\documentclass[11pt]{article}
\usepackage[utf8]{inputenc}
\usepackage[english]{babel}
\usepackage{graphicx}
\usepackage{amsmath}
\usepackage{amssymb}
\usepackage{amsfonts}
\usepackage{amsthm}
\usepackage[left=2.1cm,top=1.5cm,right=2.1cm, bottom=2.2cm,letterpaper]{geometry}
\usepackage{mathrsfs}
\usepackage{caption}
\usepackage{subcaption}
\usepackage{latexsym}
\usepackage{xfrac}
\usepackage{tcolorbox}
\usepackage{mleftright}
\usepackage{enumitem}
\usepackage{float}
\usepackage{hyperref}
\usepackage[all]{hypcap}
\usepackage{autonum}
 \usepackage{url}
\usepackage[toc]{appendix}
\newtheorem{theorem}{Theorem}[section]
\newtheorem{corollary}{Corollary}[section]
\newtheorem{lemma}{Lemma}[section]
\newtheorem{definition}{Definition}[section]

\newtheorem{remark}{Remark}[section]
\numberwithin{equation}{section}
\numberwithin{figure}{section}

\newcommand{\R}{\mathbb{R}}

\makeatletter
\renewcommand*\env@matrix[1][*\c@MaxMatrixCols c]{%
	\hskip -\arraycolsep
	\let\@ifnextchar\new@ifnextchar
	\array{#1}}
\makeatother

\def\neweq#1{\begin{equation}\label{#1}}
\def\endeq{\end{equation}}

\begin{document}
	
\title{\vspace{-0.5cm}On the steady motion of a Navier-Stokes flow across a sieve with prescribed pressure drop in a finite pipe}

\author{Gianmarco Sperone\\
{\small Faculty of Mathematics, Pontifical Catholic University of Chile}}
\date{}
\maketitle
\vspace*{-6mm}
\begin{abstract}
	\noindent
The steady motion of a viscous incompressible fluid through a sieve (that is, a wall perforated with a large number of small holes), in a pipe of finite length, is modeled through the Navier-Stokes equations under mixed boundary conditions involving the Bernoulli pressure and the tangential velocity on the inlet and outlet of the tube, while the pressure drop is prescribed along the pipe. Applying the classical energy method in homogenization theory, we study the asymptotic behavior of the solutions to this system, without any restriction on the magnitude of the data, as the diameters of the perforations vanish. Regardless of the initial scaling and distribution of the holes, we show that the sieve asymptotically becomes a wall, meaning that the effective equations are two, independent, stationary Navier-Stokes systems with a no-slip boundary condition on the wall. In the absence of external forces we prove, furthermore, that the fluid motion becomes quiescent in the homogenization limit.
\end{abstract}	

\section{Introduction and presentation of the problem}
\begin{figure}[H]
	\begin{center}
		\includegraphics[scale=0.85]{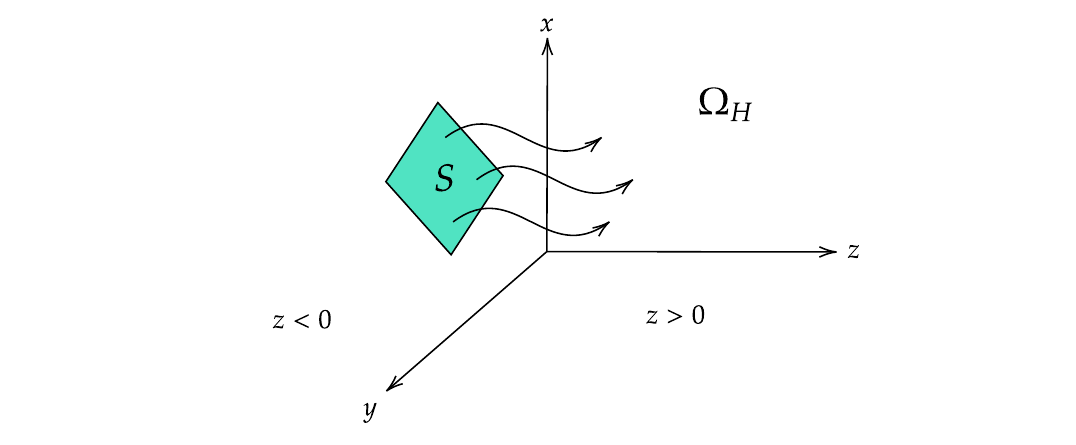}
	\end{center}
	\vspace*{-8mm}
	\caption{Representation of an \textit{aperture domain} in the sense of Heywood.}\label{dom0}
\end{figure}

In his celebrated paper from 1976, concerning uniqueness questions in the theory of viscous flow \cite{heywood1976uniqueness}, Heywood established the problem of determining the steady motion of a viscous incompressible fluid in a three-dimensional \textit{aperture domain}, that is, in an unbounded region (with a non-compact boundary) of the following type:
$$
\Omega_{H} \doteq \left\lbrace (x,y,z) \in \mathbb{R}^{3} \ \vert \ z \neq 0 \quad \text{or} \quad (x,y) \in S \, \right\rbrace \footnote{ \noindent In \cite[Section 6]{heywood1976uniqueness} Heywood showed that, for such aperture domain $\Omega_{H}$, the completion, in the Dirichlet norm, of the space of compactly supported solenoidal vector fields in $\mathcal{C}^{\infty}(\Omega_{H})$ \textbf{does not} coincide with the set of divergence-free vector fields in the Sobolev space $W_{0}^{1,2}(\Omega_{H})$.} \, ,
$$
where $S \subset \mathbb{R}^3$ is a bounded and connected subset of the $XY$-plane having a Lipschitz boundary; in words, the fluid is allowed to flow between the half-space $z<0$ and the half-space $z>0$ only through a \textit{perforation} located on the $XY$-plane, defined by the set $S$, see Figure \ref{dom0} below for a simple illustration.
\noindent
Among several other issues, Heywood discovered \cite[Section 6]{heywood1976uniqueness} (see also \cite{heywood2006auxiliary}) that the hydrodynamic equations governing the motion of the \textit{jet} can be correctly formulated by prescribing, as auxiliary conditions, either the \textit{transverse flux} through $S$ or the \textit{pressure drop} across $S$ (further properties of such steady jets can be found in the works of Borchers, Galdi \& Pileckas \cite{borchers1993uniqueness,borchers1992existence}). Later on, the use of these physical quantities as admissible boundary conditions for the Stokes or Navier-Stokes equations (in both steady and unsteady regimes) in bounded domains comprised the main subject of various investigations, of numerical and analytical nature, such as \cite{al2017instationary,concaenglish,conca1995navier,heywood1996artificial,korobkov2020solvability,maruvsic2007navier,nowakowski2020flow}. This raises a first natural question:
\begin{tcolorbox} 
\begin{center}
How does the motion of the jet change if we allow it to transit through several apertures, driven by a pressure difference across the wall containing the holes?
\end{center}
\end{tcolorbox}
\noindent
Furthermore,
\begin{tcolorbox} 
\begin{center}	
Can we describe the asymptotic motion (percolation) of the pressure-driven jet when, simultaneously, the number of apertures goes to infinity but their diameters vanish?
\end{center}
\end{tcolorbox}
\noindent
Naturally, these two inquiries should be tackled in the framework of homogenization theory concerning \textit{sieve problems}, where a partial differential equation is studied in, say, a three-dimensional domain whose intersection with a given plane contains a large number of tiny perforations. The objective is then to describe the limiting behavior of the solutions to such PDE as the sieve or grill becomes finer and finer, yielding an \textit{effective equation} in the homogeneous medium. Starting from the seminal work of Marchenko \& Khruslov \cite{marchenko1964boundary} (consult also their book \cite{marchenko1974boundary}), homogenization problems involving sieves have been analyzed predominantly in the realm of elliptic PDE's, see \cite{cioranescu2008periodic,damlamian1985probleme,del1987thick,murat1985neumann,nguetseng1985problemes,picard1987analyse,sanchez1982boundary} and references therein for a non-exhaustive list of classic works. Analogous percolation processes, in the context of Fluid Mechanics of incompressible flows, were thoroughly examined in the benchmark papers of Allaire \cite{allaire2}, Conca \cite{concasieve1,concasieve2} and Sánchez-Palencia \cite{sanchez1985probleme} for the (linear) steady-state Stokes system with non-homogeneous Dirichlet boundary conditions; further investigations in the linear setting were also conducted in \cite{bourgeat1995effective,brillard2016boundary,gomez2018homogenization,maris2012effective,maruvsic1998low}. It would appear that far less attention has been devoted to the nonlinear counterpart, and we highlight the work of J{\"a}ger \& Mikeli{\'c} \cite{jager1998effective} relating to the time-dependent Navier-Stokes equations in a bounded domain containing a \textit{thick} grill. The present article intends to serve as a contribution in the analysis of the Navier-Stokes sieve problem, with the
purpose of conferring a deeper understanding of the questions described before and motivating ulterior studies on this matter.
\par 
Let us now describe in detail the mathematical problem we are concerned with. Given $h>R>0$, in the space $\R^3$ we consider an open straight cylinder $\Omega$ of radius $R$ and length $2h$ whose axis of symmetry is directed along the $z$-axis:
$$
\Omega \doteq \left\lbrace (x,y,z) \in \mathbb{R}^{3} \ \vert \ 0 \leq x^{2} + y^2 < R^2 \, , \ -h < z < h \, \right\rbrace \, .
$$
For any $\xi \in \mathbb{R}^{2}$ and $r > 0$ we denote by $D(\xi, r) \subset \mathbb{R}^{2}$ the open disk of radius $r$ with center at $\xi$. In particular, define $\Sigma \doteq D((0,0),R)$ as the middle cross-section of $\Omega$. Let $( K_{n} )_{n \in \mathbb{N}}$ be a sequence of open, bounded and simply connected subsets of $\mathbb{R}^{2}$ with a $\mathcal{C}^{2}$-boundary such that $(0,0) \in K_{n}$, $\forall n \in \mathbb{N}$, and
$$
\sup_{n \in \mathbb{N}} | K_{n} | < \infty \, .
$$
Take $\varepsilon_{*} \in (0,1)$ such that $\varepsilon_{*} | K_{n} | < \pi R^{2}$, $\forall n \in \mathbb{N}$. Following \cite{nevcasova2023homogenization,nevcasova2022homogenization}, given $\alpha > 0$ and $\varepsilon \in (0,\varepsilon_{*}]$, define $r_{\varepsilon} \doteq \text{exp}(-\varepsilon^{-\alpha})$ and suppose there exist an integer $N(\varepsilon) \geq 1$ and a collection of points $\xi^{\varepsilon}_{1},...,\xi^{\varepsilon}_{N(\varepsilon)} \in \mathbb{R}^{2}$ satisfying the following properties:
\begin{equation} \label{perforation}
	\begin{aligned}
		& \xi^{\varepsilon}_{n} + r_{\varepsilon} \, \overline{K_{n}} \subset D(\xi^{\varepsilon}_{n}, \delta_{0} \, r_{\varepsilon}) \subset D(\xi^{\varepsilon}_{n}, \delta_{1} \, \varepsilon) \subset \Sigma \qquad \forall n \in \{1,...,N(\varepsilon)\} \,, \\[6pt]
		& \partial D \left(\xi^{\varepsilon}_{n}, \delta_{1} \, \varepsilon \right) \cap \partial D \left(\xi^{\varepsilon}_{m}, \delta_{1} \, \varepsilon \right) = \emptyset \qquad \forall n,m \in \{1,...,N(\varepsilon)\} \,, \ n \neq m \, , \\[6pt]
		& \partial D \left(\xi^{\varepsilon}_{n}, \delta_{1} \, \varepsilon \right) \cap \partial \Sigma = \emptyset \qquad \forall n \in \{1,...,N(\varepsilon)\} \, ,
	\end{aligned}
\end{equation}
for some constants $\delta_{0}, \delta_{1} > 0$ that are independent of $\varepsilon \in (0,\varepsilon_{*}]$. Setting $K^{\varepsilon}_{n} \doteq \xi^{\varepsilon}_{n} + r_{\varepsilon} \, K_{n}$ for every $n \in \{1,...,N(\varepsilon)\}$, we will refer to the family $\{ K^{\varepsilon}_{n} \}^{N(\varepsilon)}_{n=1}$ satisfying \eqref{perforation} as the \textit{perforations} or \textit{apertures}, while the set
\begin{equation} \label{sieve}
	\Gamma_{\varepsilon} \doteq \Sigma \setminus \overline{K_{\varepsilon}} \doteq \Sigma \setminus \bigcup^{N(\varepsilon)}_{n=1} \overline{K^{\varepsilon}_{n}} \, ,
\end{equation}
represents the \textit{sieve} or \textit{grill} at the $\varepsilon$-level. We emphasize that, given $\varepsilon \in (0,\varepsilon_{*}]$, the family of apertures $\{ K^{\varepsilon}_{n} \}^{N(\varepsilon)}_{n=1}$ is built in such a way that the \textit{size} of each perforation is proportional to $r_{\varepsilon}$, while the mutual distance between any two consecutive holes is proportional to $\varepsilon$. Moreover, since we only consider those apertures that are \textit{strictly} contained in $\Sigma$ (in the sense of \eqref{perforation}$_3$), the following bound on the number $N(\varepsilon)$ holds:
\begin{equation} \label{numero}
N(\varepsilon) \leq \left( \dfrac{R}{\delta_{1} \, \varepsilon} \right)^2 \, .
\end{equation}
Notice, however, that the sets $\{ K^{\varepsilon}_{n} \}^{N(\varepsilon)}_{n=1}$ may have different shapes and that they are not necessarily periodically distributed in $\Sigma$. We now define the \textit{fluid domain} (open, connected and bounded set) at the $\varepsilon$-level as
\begin{equation} \label{perfordomain}
	\Omega_{\varepsilon} \doteq \Omega \setminus \overline{\Gamma_{\varepsilon}} = \Omega_{-} \cup K_{\varepsilon} \cup \Omega_{+} \, ,
\end{equation}
where
\begin{equation} \label{perfordomainpm}
	\Omega_{-} \doteq \left\lbrace (x,y,z) \in \Omega \ \vert \ -h < z < 0 \, \right\rbrace \qquad \text{and} \qquad \Omega_{+} \doteq \left\lbrace (x,y,z) \in \Omega \ \vert \ 0 < z < h \, \right\rbrace \, .
\end{equation}
The boundary of $\Omega_{\varepsilon}$ is decomposed as $ \partial \Omega_{\varepsilon} = \Gamma_{I} \cup \Gamma^{\varepsilon}_{W} \cup \Gamma_{O}$, where
\begin{equation}\label{boundaryomega1}
	\begin{aligned}
		& \Gamma_{I} = \left\lbrace (x,y,z) \in \mathbb{R}^{3} \ \vert \  x^{2} + y^2 < R^2, \ z=-h \, \right\rbrace \, , \\[6pt]
		& \Gamma^{\varepsilon}_{W} = \mathcal{L} \cup \Gamma_{\varepsilon} \doteq \left\lbrace (x,y,z) \in \mathbb{R}^{3} \ \vert \ x^{2} + y^2 = R^2 \, , \ -h < z < h \, \right\rbrace  \cup \Gamma_{\varepsilon} \, , \\[6pt]
		& \Gamma_{O} = \left\lbrace (x,y,z) \in \mathbb{R}^{3} \ \vert \  x^{2} + y^2 < R^2, \ z=h \, \right\rbrace \, .
	\end{aligned}
\end{equation}
The outward unit normal to $\partial \Omega_{\varepsilon}$ is denoted by $\nu$ (with some abuse of notation, as such vector also depends on $\varepsilon$). Henceforth we will refer to $\Gamma_{I}$ and $\Gamma_{O}$ in \eqref{boundaryomega1} as the \textit{inlet} and \textit{outlet} of $\Omega$, respectively, while $\Gamma^{\varepsilon}_{W}$ includes all the \textit{physical walls} of $\Omega_{\varepsilon}$, see Figure \ref{dom1} below.

\begin{figure}[H]
	\begin{center}
		\includegraphics[height=45mm,width=170mm]{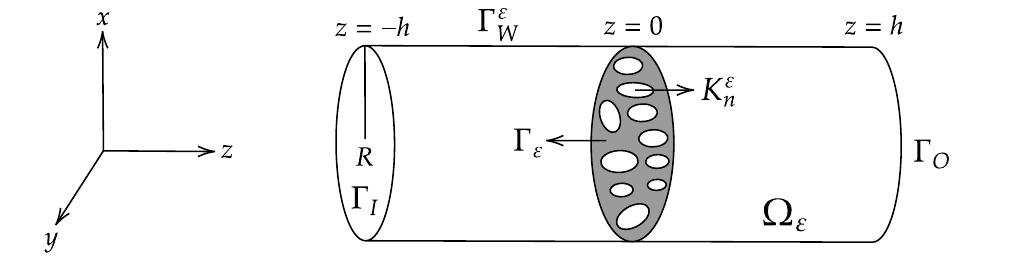}
	\end{center}
	\vspace*{-5mm}
	\caption{Representation of the fluid domain $\Omega_{\varepsilon}$ and the sieve $\Gamma_{\varepsilon}$.}\label{dom1}
\end{figure}
\noindent
Given $\varepsilon \in (0,\varepsilon_{*}]$, we analyze the steady motion of a viscous incompressible fluid (having a constant kinematic viscosity equal to 1) along $\Omega_{\varepsilon}$, which is characterized by its velocity vector field $u_{\varepsilon} : \Omega_{\varepsilon} \longrightarrow \mathbb{R}^3$ and its scalar pressure $p_{\varepsilon} : \Omega_{\varepsilon} \longrightarrow \mathbb{R}$, under the action of an external force $f : \Omega \longrightarrow \mathbb{R}^3$. Such stationary motion will be modeled through the \textbf{prescribed pressure drop problem} in $\Omega_{\varepsilon}$, that is, the following boundary-value problem (with mixed boundary conditions) associated to the steady-state Navier-Stokes equations in $\Omega_{\varepsilon}$:
\begin{equation}\label{nsstokespd}
	\left\{
	\begin{aligned}
		& -\Delta u_{\varepsilon}+(u_{\varepsilon}\cdot\nabla)u_{\varepsilon}+\nabla p_{\varepsilon}=f \, , \quad  \nabla\cdot u_{\varepsilon}=0 \ \ \mbox{ in } \ \ \Omega_{\varepsilon} \, , \\[4pt]
		& u_{\varepsilon}=0 \ \ \mbox{ on } \ \ \Gamma^{\varepsilon}_{W} \, , \\[4pt]
		& u_{\varepsilon} \times \nu = 0 \, , \quad p_{\varepsilon} +  \dfrac{1}{2} |u_{\varepsilon}|^{2} = p^{-} \ \ \mbox{ on } \ \ \Gamma_{I} \, , \\[4pt]
		& u_{\varepsilon} \times \nu = 0 \, , \quad p_{\varepsilon} +  \dfrac{1}{2} |u_{\varepsilon}|^{2} = p^{+} \ \ \mbox{ on } \ \ \Gamma_{O} \, .
	\end{aligned}
	\right.
\end{equation}
While \eqref{nsstokespd}$_{2}$ describes the usual no-slip boundary condition on the physical walls $\Gamma^\varepsilon_{W}$, the first equality in \eqref{nsstokespd}$_{3}$-\eqref{nsstokespd}$_{4}$ dictates that the fluid must enter and leave the domain $\Omega$ orthogonal to the inlet and outlet walls. The second identity in \eqref{nsstokespd}$_{3}$-\eqref{nsstokespd}$_{4}$ imposes that, respectively on the inlet $\Gamma_{I}$ and outlet $\Gamma_{O}$, the \textit{Bernoulli pressure} defined as $\Phi_{\varepsilon} \doteq p_{\varepsilon} + |u_{\varepsilon}|^{2}/2$ must equal some prescribed constants $p^{\mp} \in \mathbb{R}$ that represent the \textit{pressure drop} $p^{+} - p^{-}$ along the pipe containing the sieve. 
Given a velocity field $u_{\varepsilon} \in \mathcal{C}^{2}(\Omega_{\varepsilon}) \cap \mathcal{C}(\overline{\Omega_{\varepsilon}})$ solving \eqref{nsstokespd}, observe that the transversal flux rate
\begin{equation} \label{refer00}
	F_{\varepsilon} \doteq \int_{\Sigma(s)} u_{\varepsilon} \cdot \widehat{k} \qquad \forall s \in [-h,h] \, ,
\end{equation}
is constant along the pipe, but \textit{depends} on the solution, where  $\Sigma(0) \doteq \Sigma$ and
\begin{equation} \label{seccionest}
\Sigma(s) \doteq \left\lbrace (x,y,z) \in \mathbb{R}^{3} \ \vert \ 0 \leq x^{2} + y^2 < R^2, \ z=s \, \right\rbrace \qquad \forall s \in [-h,h] \setminus \{0\} \, .
\end{equation}
Indeed, given $s \in (0,h]$, we define the region
$$
\Omega_{\varepsilon}(s) \doteq \left\lbrace (x,y,z) \in \Omega_{\varepsilon} \ \vert -h < z < s \, \right\rbrace \, .
$$
Since $u_{\varepsilon}$ vanishes on $\Gamma^{\varepsilon}_{W}$, after applying the Divergence Theorem (separately in $\Omega_{-}$ and $\Omega_{\varepsilon}(s) \cap \Omega_{+}$) we infer the following:
$$
\begin{aligned}
0 & = \int_{\Omega_{\varepsilon}(s)} \nabla \cdot u_{\varepsilon} = \int_{\Omega_{-}} \nabla \cdot u_{\varepsilon} + \int_{\Omega_{\varepsilon}(s) \cap \Omega_{+}} \nabla \cdot u_{\varepsilon} \\[6pt]
& = -\int_{\Gamma_{I}} u_{\varepsilon} \cdot \widehat{k} + \int_{K_{\varepsilon}} u_{\varepsilon} \cdot \widehat{k} + \int_{\Sigma(s)} u_{\varepsilon} \cdot \widehat{k} - \int_{K_{\varepsilon}} u_{\varepsilon} \cdot \widehat{k} = \int_{\Sigma(s)} u_{\varepsilon} \cdot \widehat{k} - \int_{\Gamma_{I}} u_{\varepsilon} \cdot \widehat{k} \, ,
\end{aligned}
$$
so that
$$
\int_{\Sigma(s)} u_{\varepsilon} \cdot \widehat{k} = \int_{\Gamma_{I}} u_{\varepsilon} \cdot \widehat{k} \qquad \forall s \in (0,h] \, .
$$
In a similar way one can prove \eqref{refer00} for every $s \in [-h,0]$.

In view of the identity
$$
\nabla \left( \dfrac{1}{2} |u_{\varepsilon}|^2 \right) = (\nabla u_{\varepsilon})^{\top} u_{\varepsilon} \quad \text{in} \quad \Omega_{\varepsilon} \, ,
$$
it is customary (see, for example, \cite{heywood1996artificial,korobkov2020solvability}) to add the term $(\nabla u_{\varepsilon})^{\top} u_{\varepsilon}$ to both sides of the equation of conservation of momentum \eqref{nsstokespd}$_1$, thereby resulting in the problem
\begin{equation}\label{nsstokespd1}
	\left\{
	\begin{aligned}
		& -\Delta u_{\varepsilon}+ (u_{\varepsilon}\cdot\nabla)u_{\varepsilon} - (\nabla u_{\varepsilon})^{\top} u_{\varepsilon} + \nabla \Phi_{\varepsilon}=f \, ,\ \quad  \nabla\cdot u_{\varepsilon}=0 \ \ \mbox{ in } \ \ \Omega_{\varepsilon} \, , \\[3pt]
		& u_{\varepsilon}=0 \ \ \mbox{ on } \ \ \Gamma^{\varepsilon}_{W} \, , \\[3pt]
		& u_{\varepsilon} \times \nu = 0 \, , \quad \Phi_{\varepsilon} = p^{-} \ \ \mbox{ on } \ \ \Gamma_{I} \, ,\\[3pt]
		& u_{\varepsilon} \times \nu = 0 \, , \quad \Phi_{\varepsilon} = p^{+} \ \ \mbox{ on } \ \ \Gamma_{O} \, .
	\end{aligned}
	\right.
\end{equation}

As already announced, the main goal of the present manuscript is to study the asymptotic behavior of the solutions of \eqref{nsstokespd1} as $\varepsilon \to 0^{+}$. Our main homogenization result reads:

\begin{theorem} \label{maintheoremintro}
	Let $(\Omega_{\varepsilon})_{\varepsilon \in I_{*}}$ be the family of domains \eqref{perfordomain} and $\alpha >0$. Given any $p^{\pm} \in \mathbb{R}$ and $f \in L^{2}(\Omega)$, let $(u_{\varepsilon},\Phi_{\varepsilon}) \in H^{1}(\Omega_{\varepsilon}) \times L^{2}(\Omega_{\varepsilon})$ be a weak solution of problem \eqref{nsstokespd1}. Then, up to the extraction of a subsequence, the sequences of restrictions $\{(u^{\pm}_{\varepsilon}, \Phi^{\pm}_{\varepsilon}) \}_{\varepsilon \in I_{*}} \subset H^{1}(\Omega_{\pm}) \times L^{2}(\Omega_{\pm})$ converge strongly to weak solutions $(u^{\pm}, \Phi^{\pm}) \in H^{1}(\Omega_{\pm}) \times L^{2}(\Omega_{\pm})$ of the following boundary-value problems, in $\Omega_{\pm}$, as $\varepsilon \to 0^{+}$:
	\begin{equation}\label{nsstokespdfinal-int}
		\left\{
		\begin{aligned}
			& -\Delta u^{-} + (u^{-}\cdot\nabla)u^{-} - (\nabla u^{-})^{\top} u^{-} + \nabla \Phi^{-} = f^{-} \, ,\ \quad  \nabla\cdot u^{-} = 0 \ \ \mbox{ in } \ \ \Omega_{-} \, , \\[3pt]
			& u^{-} \times \nu = 0 \, , \quad \Phi^{-} = p^{-} \ \ \mbox{ on } \ \ \Gamma_{I} \, ,\\[3pt]
			& u^{-}=0 \ \ \mbox{ on } \ \ (\partial \Omega_{-} \cap \mathcal{L}) \cup \Sigma \, ,
		\end{aligned}
		\right.
	\end{equation}
	and
	\begin{equation}\label{nsstokespdfinal+int}
		\left\{
		\begin{aligned}
			& -\Delta u^{+} + (u^{+}\cdot\nabla)u^{+} - (\nabla u^{+})^{\top} u^{+} + \nabla \Phi^{+} = f^{+} \, ,\ \quad  \nabla\cdot u^{+} = 0 \ \ \mbox{ in } \ \ \Omega_{+} \, , \\[3pt]
			& u^{+} \times \nu = 0 \, , \quad \Phi^{+} = p^{+} \ \ \mbox{ on } \ \ \Gamma_{O} \, ,\\[3pt]
			& u^{+}=0 \ \ \mbox{ on } \ \ (\partial \Omega_{+} \cap \mathcal{L}) \cup \Sigma \, .
		\end{aligned}
		\right.
	\end{equation}	
\end{theorem} 
\noindent
In order to prove Theorem \ref{maintheoremintro}, we firstly ensure the unrestricted weak solvability of problem \eqref{nsstokespd1} (concerning the size of the data of the problem, that is, external force and prescribed pressure drop) for a fixed $\varepsilon \in (0,\varepsilon_{*}]$. This is achieved in Theorem \ref{epslevel} by means of the Leray-Schauder Principle, following the outlines of \cite[Theorem 3.2]{korobkov2020solvability}. Alongside this result, in Section \ref{epslevelsec} we provide a series of lemmas that allow us, in turn, to derive the uniform $\varepsilon$-independent bounds for the solutions of \eqref{nsstokespd1} (Theorem \ref{epslevelbounds}) that are necessary to produce the homogenization limit in Theorem \ref{maintheoremintro}. Among those lemmas, we underline the outcome of Lemma \ref{bogeps} guaranteeing the unbounded growth, as $\varepsilon \to 0^{+}$, of the continuity constant of the right inverse of the divergence operator in $\Omega_{\varepsilon}$. One of the main difficulties of our problem arises as a consequence of this, since such continuity constant is essential for the obtainment of pressure estimates; we are then led to separate our analysis into each of the sub-domains $\Omega_{\pm}$, determining the projection of the Bernoulli pressure $\Phi_{\varepsilon}$ onto the spaces of square-integrable functions having zero mean value in $\Omega_{\pm}$. Subsequently, applying the classical energy method in homogenization theory \cite[Appendix]{sanchez1980non}, in Section \ref{energymethod} (specifically, Theorem \ref{effectiveq1}) we show that, as $\varepsilon \to 0^{+}$, the sieve $\Gamma_{\varepsilon}$ asymptotically becomes the \textit{physical} wall $\Sigma$, meaning that the effective equations are the two, independent, stationary Navier-Stokes systems \eqref{nsstokespdfinal-int}-\eqref{nsstokespdfinal+int} with a no-slip boundary condition on $\Sigma$. Quite surprisingly, this result holds for any value of $\alpha > 0$, that is, regardless of the initial size of the perforations. This is in clear contrast with other existing results in the homogenization literature, where three asymptotic outcomes can be obtained depending upon the value of $\alpha > 0$, see \cite{allaire2,concasieve1,fernandez2008numerical,gomez2018homogenization,maris2012effective,picard1987analyse,sanchez1982acoustic}. Roughly speaking, a priori one would expect to generate the following effective laws \footnote{ \noindent The determination of $\alpha =2$ as the critical value follows from the estimates, with respect to $\varepsilon \in (0,\varepsilon_{*}]$, concerning the two-dimensional \textit{relative capacity} of the perforations $K_{\varepsilon}$ inside the disk $\Sigma$, as in \cite[Lemma 2.1]{patriarca2023homogenization}.}:
\begin{itemize}[leftmargin=3mm]
	\item[-] In the sub-critical case, when $0 < \alpha < 2$, the perforations are too large, and therefore the presence of the sieve does not affect the fluid motion. The homogenized equations then coincide with those of the original prescribed pressure drop problem in the whole pipe $\Omega$.
	\item[-] In the critical case, when $\alpha = 2$, the homogenized equations exhibit an \textit{unexpected} non-homogeneous Dirichlet boundary condition on $\Sigma$ (leakage), see \cite[Section 2.2]{concasieve1}.
	\item[-] In the super-critical case, when $\alpha > 2$, the perforations are too small, so that the sieve \textit{converges} to the wall $\Sigma$ and the homogenized equations concur precisely with \eqref{nsstokespdfinal-int}-\eqref{nsstokespdfinal+int}.
\end{itemize}
Nevertheless, special attention should be drawn to \cite[Théorème 2.3 bis]{concasieve1}, where the Author proves a result analogous to Theorem \ref{maintheoremintro} in the critical case when $\alpha = 2$, but assuming that the tranverse flux across $\Sigma$ (which, under non-homogeneous Dirichlet boundary conditions, is independent of $\varepsilon \in (0,\varepsilon_{*}]$) is null. In our setting, the $\varepsilon$-uniform bounds obtained in Theorem \ref{epslevelbounds}, combined with an $\varepsilon$-dependent trace inequality on $\Sigma$ satisfied by every weak solution of \eqref{nsstokespd1}, enable us to show that $F_{\varepsilon}$ decays as $\sqrt{r_{\varepsilon}}$ as $\varepsilon \to 0^{+}$; see Lemma \ref{concalem} and Corollary \ref{concacor}.

\newpage
\section{Solvability of the boundary-value problem at the $\varepsilon$-level} \label{epslevelsec}
Let $\varepsilon \in I_{*}$ be a fixed parameter, with $I_{*} \doteq (0,\varepsilon_{*}]$. As in \cite[Section 2]{sperone2023}, we introduce the functional spaces (of vector fields) that will be employed hereafter:
$$
\mathcal{V}(\Omega_{\varepsilon}) \doteq \left\lbrace v \in H^{1}(\Omega_{\varepsilon}) \ | \ v \times \nu = 0 \ \ \mbox{on} \ \ \Gamma_{I} \cup \Gamma_{O} \, ; \qquad v = 0 \ \ \mbox{on} \ \ \Gamma^{\varepsilon}_{W} \, \right\rbrace
$$
and
$$
\mathcal{V}_{\sigma}(\Omega_{\varepsilon}) \doteq \left\lbrace v \in H^{1}(\Omega_{\varepsilon}) \ | \ \nabla \cdot v=0 \ \ \mbox{in} \ \ \Omega_{\varepsilon} \, ; \qquad v \times \nu = 0 \ \ \mbox{on} \ \ \Gamma_{I} \cup \Gamma_{O} \, ; \qquad v = 0 \ \ \mbox{on} \ \ \Gamma^{\varepsilon}_{W} \, \right\rbrace \, ,
$$
which are Hilbert spaces if endowed with the Dirichlet scalar product of the gradients, denoted by
\begin{equation} \label{dirscal}
	[v, w]_{\mathcal{V}(\Omega_{\varepsilon})} \doteq \int_{\Omega_{\varepsilon}} \nabla v \cdot \nabla w \qquad \forall v,w \in \mathcal{V}(\Omega_{\varepsilon}) \, .
\end{equation}

\begin{remark} \label{nonlip}
	Notice that $\Omega_{\varepsilon}$ is \textbf{not} a locally Lipschitz domain: in a neighborhood of the perforated wall $\Gamma_{\varepsilon}$, $\Omega_{\varepsilon}$ is not locally placed on one side of its boundary. However, the domains $\Omega_{\pm}$ are Lipschitz, see \eqref{perfordomainpm}, so that the elements of $H^{1}(\Omega_{\pm})$ have a well-defined trace on $\Sigma$. Accordingly, to fix notation and following \cite{concasieve1,concasieve2}, given any function $v \in L^{1}(\Omega_{\varepsilon})$ (scalar or vector), we denote by $v^{\pm} \in L^{1}(\Omega_{\pm})$ its restriction to $\Omega_{\pm}$. It can then be easily proved (as in \cite[Theorem 1.7.1]{grisvard2011elliptic}) that, for any $v \in L^{2}(\Omega_{\varepsilon})$ (scalar or vector), there holds:
	$$
	v \in H^{1}(\Omega_{\varepsilon}) \qquad \Longleftrightarrow \qquad v^{\pm} \in H^{1}(\Omega_{\pm}) \quad \text{and} \quad v^{+} |_{K_{\varepsilon}} = v^{-} |_{K_{\varepsilon}} \, .
	$$
	In particular, for any $v \in L^{2}(\Omega_{\varepsilon})$ (scalar or vector), it follows that
	$$
	v \in H_{0}^{1}(\Omega_{\varepsilon}) \qquad \Longleftrightarrow \qquad v^{\pm} \in H^{1}(\Omega_{\pm}) \, , \quad v^{+} |_{K_{\varepsilon}} = v^{-} |_{K_{\varepsilon}} \quad \text{and} \quad v |_{\partial \Omega} = v^{\pm} |_{\Gamma_{\varepsilon}} = 0 \, .
	$$
	Nevertheless, given $v \in H^{1}(\Omega_{\varepsilon})$ (scalar or vector), the traces $v^{\pm} |_{\Gamma_{\varepsilon}} \in H^{1/2}(\Gamma_{\varepsilon})$ do not necessarily coincide almost everywhere on $\Gamma_{\varepsilon}$.
\end{remark}

Before proving the weak solvability of problem \eqref{nsstokespd1} (equivalently, of problem \eqref{nsstokespd}), we will state a series of lemmas that shall play a key role in the sequel. The first (elementary) one concerns the uniform boundedness (with respect to $\varepsilon \in I_{*}$) of some Sobolev embedding constants:

\begin{lemma} \label{unisob}
	Let $\Omega_{\varepsilon}$ be as in \eqref{perfordomain}. For any $p \in [1,6]$ and any function (scalar or vector) $\varphi \in H^{1}(\Omega_{\varepsilon})$ that vanishes on $\mathcal{L}$ there holds the estimate
	\begin{equation} \label{sob0}
		\| \varphi \|_{L^{p}(\Omega_{\varepsilon})} \leq C_{*} \| \nabla \varphi \|_{L^{2}(\Omega_{\varepsilon})} \qquad \forall \varepsilon \in I_{*} \, ,
	\end{equation}
	for some constant $C_{*} > 0$ that depends on $\Omega$, $p \in [1,6]$ and $\{ \delta_{0}, \delta_{1} \}$, but is independent of $\varepsilon \in I_{*}$.	
\end{lemma}
\noindent
\begin{proof}
Given $p \in [1,6]$, in what follows, $C > 0$ will always denote a generic constant that depends on $\Omega$, $p \in [1,6]$ and $\{ \delta_{0}, \delta_{1} \}$ (independently of $\varepsilon \in I_{*}$), but that may change from line to line.
\par
Given any function (scalar or vector) $\varphi \in H^{1}(\Omega_{\varepsilon})$ vanishing on $\mathcal{L}$, since $\varphi^{\pm} \in H^{1}(\Omega_{\pm})$ vanishes on $\mathcal{L} \cap \partial \Omega_{\pm}$, we can apply the Sobolev inequality, separately in $\Omega_{\pm}$, to deduce that
$$
\|\varphi\|^{p}_{L^{p}(\Omega_{\varepsilon})} = \|\varphi^{-}\|^{p}_{L^{p}(\Omega_{-})} + \|\varphi^{+}\|^{p}_{L^{p}(\Omega_{+})} \leq C \left( \| \nabla \varphi^{-}\|^{p}_{L^{2}(\Omega_{-})} + \| \nabla \varphi^{+}\|^{p}_{L^{2}(\Omega_{+})} \right) \leq 2 C \| \nabla \varphi\|^{p}_{L^{2}(\Omega_{\varepsilon})} \, ,
$$
which finishes the proof.
\end{proof}

Our second lemma in this section, in the spirit of \cite[Lemme 3.2]{concasieve1}, involves a $\varepsilon$-dependent trace inequality for functions that vanish on the wall $\Gamma_{\varepsilon}$.

\begin{lemma} \label{concalem}
	Let $\Omega_{\varepsilon}$ be as in \eqref{perfordomain}. For any function (scalar or vector) $\varphi \in H^{1}(\Omega_{\varepsilon})$ that vanishes on $\Gamma_{\varepsilon}$ there holds the estimate
	\begin{equation} \label{conca0}
		\| \varphi \|_{L^{2}(\Sigma)} \leq C_{*} \sqrt{r_{\varepsilon}} \, \| \nabla \varphi \|_{L^{2}(\Omega_{\varepsilon})} \qquad \forall \varepsilon \in I_{*} \, ,
	\end{equation}
	for some constant $C_{*} > 0$ that depends on $\Omega$ and $\{ \delta_{0}, \delta_{1} \}$, but is independent of $\varepsilon \in I_{*}$.	
\end{lemma}
\noindent
\begin{proof}
	In what follows, $C > 0$ will always denote a generic constant that depends on $\Omega$ and $\{ \delta_{0}, \delta_{1} \}$ (independently of $\varepsilon \in I_{*}$), but that may change from line to line.
	\par
Given any $n \in \{1,...,N(\varepsilon)\}$, let $\{ (a^{\varepsilon}_{n}, b^{\varepsilon}_{n}) \}^{N(\varepsilon)}_{n=1} \subset \mathbb{R}^2$ be such that $\xi^{\varepsilon}_{n} = (a^{\varepsilon}_{n}, b^{\varepsilon}_{n})$, see again \eqref{perforation}. Define then the (open) cylinders
$$
\mathcal{P}^{\varepsilon}_{n} \doteq \left\lbrace (x,y,z) \in \mathbb{R}^{3} \ \Bigg\vert \ \left( \dfrac{x - a^{\varepsilon}_{n}}{r_{\varepsilon}} \right)^{2} + \left( \dfrac{y - b^{\varepsilon}_{n}}{r_{\varepsilon}} \right)^{2} < \delta_{0}^{2} \, , \ 0 < z < r_{\varepsilon} \, h \, \right\rbrace \qquad \forall n \in \{1,...,N(\varepsilon)\} \, ,
$$
and also $\mathcal{P}_{0} \doteq D((0,0), \delta_{0}) \times (0,h)$. Therefore, given any (scalar or vector) $\varphi \in H^{1}(\Omega_{\varepsilon})$ vanishing on $\Gamma_{\varepsilon}$, the trace inequality 
\begin{equation} \label{thintube1}
\| \varphi \|_{L^{2}(D((0,0), \delta_{0}))} \leq C \| \nabla \varphi^{+} \|_{L^{2}(\mathcal{P}_{0} )} \, 
\end{equation}
certainly holds, see Remark \ref{nonlip}. Applying \eqref{thintube1} and the changes of variables
$$
\begin{aligned}
(x,y) \in D \left(\xi^{\varepsilon}_{n}, \delta_{0} \, r_{\varepsilon} \right) & \longleftrightarrow \left( \dfrac{x - a^{\varepsilon}_{n}}{r_{\varepsilon}} , \dfrac{y - b^{\varepsilon}_{n}}{r_{\varepsilon}} \right) \in D((0,0), \delta_{0}) \, , \\[6pt]
(x,y,z) \in \mathcal{P}^{\varepsilon}_{n} & \longleftrightarrow \left( \dfrac{x - a^{\varepsilon}_{n}}{r_{\varepsilon}} , \dfrac{y - b^{\varepsilon}_{n}}{r_{\varepsilon}} , \dfrac{z}{r_{\varepsilon}} \right) \in \mathcal{P}_{0} \, ,
\end{aligned}
$$
we obtain
$$
\begin{aligned}
\| \varphi \|^{2}_{L^{2}(\Sigma)} & = \sum_{n=1}^{N(\varepsilon)} \left( \int_{K^{\varepsilon}_{n}} | \varphi(x,y,0) |^{2} \, dx \, dy \right) \leq \sum_{n=1}^{N(\varepsilon)} \left( \int_{D \left(\xi^{\varepsilon}_{n}, \delta_{0} \, r_{\varepsilon} \right)} | \varphi(x,y,0) |^{2} \, dx \, dy \right) \\[6pt]
& = r^{2}_{\varepsilon} \sum_{n=1}^{N(\varepsilon)} \left( \int_{D \left((0,0), \delta_{0} \right)} | \varphi( r_{\varepsilon} \, x + a^{\varepsilon}_{n}, r_{\varepsilon} \, y + b^{\varepsilon}_{n} , 0) |^{2} \, dx \, dy \right) \\[6pt]
& \leq C \, r^{4}_{\varepsilon} \sum_{n=1}^{N(\varepsilon)} \left( \int_{\mathcal{P}_{0}} | \nabla \varphi^{+}( r_{\varepsilon} \, x + a^{\varepsilon}_{n}, r_{\varepsilon} \, y + b^{\varepsilon}_{n} , r_{\varepsilon} \, z) |^{2} \, dx \, dy \, dz \right) \\[6pt]
& = C \, r_{\varepsilon} \sum_{n=1}^{N(\varepsilon)} \left( \int_{\mathcal{P}^{\varepsilon}_{n}} | \nabla \varphi^{+}( x, y, z) |^{2} \, dx \, dy \, dz \right) \leq C \, r_{\varepsilon} \| \nabla \varphi \|^{2}_{L^{2}(\Omega_{\varepsilon})} \, ,
\end{aligned}
$$
thus concluding the proof.
\end{proof}

Let $Q \subset \mathbb{R}^{3}$ be any bounded domain that can be decomposed as the union of a finite number of open sets, each one being star-shaped with respect to some ball strictly contained in it \footnote{This condition is satisfied, for example, by any bounded domain in $\mathbb{R}^{3}$ enjoying the \textit{interior cone property}, consult \cite[Section 2.3.4]{necas2011direct} and also \cite[Remark III.3.4]{galdi2011introduction}.}. Consider the space of square-integrable functions in $Q$ having zero mean value:
\begin{equation} \label{l02}
	L^2_0(Q) \doteq \left\{ g\in L^2(Q) \ \Big | \ \int_Q g = 0 \right\} \, .
\end{equation}
We define the \textit{Bogovskii constant} of $Q$ as
\begin{equation}\label{bogo} 
	C_B (Q) \doteq   \sup_ {g \in L ^ 2 _0 (Q) \setminus \{ 0 \} }  \inf  \left\{ \frac{  \|\nabla v\|_{L^2(Q)} }{ \|g\|_{L^2(Q)}} \ \Bigg| \ v \in H ^ 1 _ 0 (Q)  \, , \ \nabla\cdot v=g \ \text{ in } \ Q \right\}  \, ,
\end{equation}
see \cite[Section 2]{gazspefra} and also \cite[Remark III.3.4]{galdi2011introduction}. Bogovskii \cite{bogovskii1979solution} showed that, given any $q \in L^2_0(Q)$, there exists a vector field $X \in H^1_0(Q)$ such that $\nabla \cdot X=q$ in $Q$ and
$$
\| \nabla X \|_{L^{2}(Q)} \leq C_{B}(Q) \| q \|_{L^{2}(Q)} \, .
$$
For a fixed $\varepsilon \in I_{*}$, the fluid domain $\Omega_{\varepsilon}$ in \eqref{perfordomain} can be certainly decomposed as the union of a finite number of open sets, each one being star-shaped with respect to some ball strictly contained in it. Our next result shows the unbounded growth of the Bogovskii constant of $\Omega_{\varepsilon}$ with respect to $\varepsilon \in I_{*}$. 
\begin{lemma} \label{bogeps}
	Let $\Omega_{\varepsilon}$ be as in \eqref{perfordomain}. There holds the lower bound
	\begin{equation} \label{bogepsineq}
		C_{B}(\Omega_{\varepsilon}) \geq \dfrac{C_{*}}{\sqrt{r_{\varepsilon}}} \qquad \forall \varepsilon \in I_{*} \, ,
	\end{equation}
	for some constant $C_{*} > 0$ that depends on $\Omega$ and $\{ \delta_{0}, \delta_{1} \}$, but is independent of $\varepsilon \in I_{*}$.	
\end{lemma}
\noindent
\begin{proof}
	In what follows, $C > 0$ will always denote a generic constant that depends on $\Omega$ and $\{ \delta_{0}, \delta_{1} \}$ (independently of $\varepsilon \in I_{*}$), but that may change from line to line.
	\par
Define the function $g_{*} \in L^2_0(\Omega)$ (see \eqref{l02}) by
$$
g_{*}(x,y,z) =
\left\{ 
\begin{array}{ll}
	-1 & \text{if} \ \ (x,y,z) \in \Omega_{-} \\[4pt]
	0 & \text{if} \ \ (x,y,z) \in \Sigma\\[4pt] 
	1 & \text{if} \ \ (x,y,z) \in \Omega_{+} \, ,
\end{array}
\right.
$$
and therefore, we also have $g_{*} \in L^2_0(\Omega_{\varepsilon})$ with $\|g_{*}\|_{L^2(\Omega_{\varepsilon})} = \sqrt{| \Omega |}$. It then suffices to prove that
$$
\|\nabla J \|_{L^2(\Omega_{\varepsilon})} \geq \dfrac{C}{\sqrt{r_{\varepsilon}}} \quad \text{for every vector field} \ J \in H^{1}_{0}(\Omega_{\varepsilon}) \ \ \text{such that} \ \ \nabla \cdot J = g_{*} \ \ \text{in} \ \ \Omega_{\varepsilon} \, .
$$
Let $J \in H^{1}_{0}(\Omega_{\varepsilon})$ be any vector field such that $\nabla \cdot J = g_{*}$ almost everywhere in $\Omega_{\varepsilon}$. In view of the Divergence Theorem we then have
$$
\int_{K_{\varepsilon}} J \cdot \widehat{k} = \int_{\Sigma} J \cdot \widehat{k} = \int_{\Omega_{-}} \nabla \cdot J^{-} = \int_{\Omega_{-}} g^{-}_{*} = - \pi R^{2} h \, .
$$
Then, applying H\"older's inequality and \eqref{conca0} we deduce
$$
\pi R^{2} h = \left| \int_{\Sigma} J \cdot \widehat{k} \right| \leq \| J \|_{L^1(\Sigma)} \leq \sqrt{\pi} \, R \, \| J \|_{L^2(\Sigma)} \leq C \sqrt{r_{\varepsilon}} \, \| \nabla J \|_{L^{2}(\Omega_{\varepsilon})} \, ,
$$
which finishes the proof.
\end{proof}

Another essential preliminary result concerns the construction of a \textit{Bogovskii-type} operator on the space $L^{2}(\Omega_{\varepsilon})$, which exploits the corresponding operator on $L_{0}^{2}(\Omega_{\varepsilon})$. Inspired by \cite[Lemma 4.2]{korobkov2020solvability} and \cite[Lemma 2.2]{sperone2023}, we have:
\begin{lemma} \label{bogtype}
	Let $\Omega_{\varepsilon}$ be as in \eqref{perfordomain} and $q \in L^{2}(\Omega_{\varepsilon})$. There exists a vector field $Y_{\varepsilon} \in \mathcal{V}(\Omega_{\varepsilon})$ such that
	\begin{equation} \label{vecje}
     \nabla \cdot Y_{\varepsilon} = q \ \ \mbox{in} \ \ \Omega_{\varepsilon} \qquad \text{and} \qquad \| \nabla Y_{\varepsilon} \|_{L^{2}(\Omega_{\varepsilon})} \leq C_{*} (1 + C_{B}(\Omega_{\varepsilon})) \| q \|_{L^{2}(\Omega_{\varepsilon})} \qquad \forall \varepsilon \in I_{*} \, ,
	\end{equation}
	for some constant $C_{*} > 0$ that depends on $\Omega$ and $\{ \delta_{0}, \delta_{1} \}$, but is independent of $\varepsilon \in I_{*}$.
\end{lemma}
\noindent
\begin{proof}
	In what follows, $C > 0$ will always denote a generic constant that depends on $\Omega$ and $\{ \delta_{0}, \delta_{1} \}$ (independently of $\varepsilon \in I_{*}$), but that may change from line to line.
	
	Consider a Hagen-Poiseuille flow having unit flow rate in $\Omega$, that is, 
	\begin{equation} \label{pouseuille}
		U_{0}(x,y,z) \doteq \dfrac{2}{\pi R^{4}} (R^2 - x^2 - x^2) \widehat{k} \qquad \forall (x,y,z) \in \overline{\Omega} \, ,
	\end{equation}
	see \cite[Chapter II]{landau} for further details. Clearly $U_{0} \in \mathcal{C}^{\infty}(\overline{\Omega})$ is divergence-free, it vanishes on $\mathcal{L}$, and $U_{0} \times \nu =0$ on $\Gamma_{I} \cup \Gamma_{O}$. Moreover,
	\begin{equation} \label{hagen1}
		\int_{\Gamma_{O}} U_{0} \cdot \widehat{k} = 1 \, .
	\end{equation}
    We then define the vector field
	\begin{equation} \label{bogtypevec}
	Q_{\varepsilon}(x,y,z) \doteq \left( \int_{\Omega_{\varepsilon}} q \right) \dfrac{z(z+h)}{2h^2} \, U_{0}(x,y,z) \qquad \forall (x,y,z) \in \Omega_{\varepsilon} \, ,
\end{equation}
which is clearly an element of $\mathcal{V}(\Omega_{\varepsilon}) \cap \mathcal{C}^{\infty}(\overline{\Omega})$ such that $Q_{\varepsilon} = 0$ on $\Gamma_{I} \cup \Sigma$ and
	\begin{equation} \label{pointwise2}
		\| \nabla Q_{\varepsilon} \|_{L^{2}(\Omega_{\varepsilon})} \leq C \| q \|_{L^{2}(\Omega_{\varepsilon})} \, .
	\end{equation}
	On the other hand, from the Divergence Theorem and \eqref{hagen1} we obtain
	$$
	\begin{aligned}
	\int_{\Omega_{\varepsilon}} \nabla \cdot Q_{\varepsilon} = \int_{\Omega_{-}} \nabla \cdot Q^{-}_{\varepsilon} + \int_{\Omega_{+}} \nabla \cdot Q^{+}_{\varepsilon} = \int_{\Gamma_{O}} Q^{+}_{\varepsilon} \cdot \widehat{k} = \int_{\Omega_{\varepsilon}} q \, ,
	\end{aligned}
	$$
	so that $q - \nabla \cdot Q_{\varepsilon} \in L^{2}_{0}(\Omega_{\varepsilon})$, see \eqref{l02}. Then, there exists a vector field $X_{\varepsilon} \in H^{1}_{0}(\Omega_{\varepsilon})$ such that
	\begin{equation} 
		\nabla \cdot X_{\varepsilon} = q - \nabla \cdot Q_{\varepsilon} \quad \text{in} \quad \Omega_{\varepsilon} \qquad \text{and} \qquad  \| \nabla X_{\varepsilon} \|_{L^{2}(\Omega_{\varepsilon})} \leq C_{B}(\Omega_{\varepsilon}) \| q - \nabla \cdot Q_{\varepsilon} \|_{L^{2}(\Omega_{\varepsilon})}  \, .
	\end{equation}
	We set $Y_{\varepsilon} \doteq X_{\varepsilon} + Q_{\varepsilon}$ which, in view of \eqref{pointwise2}, is an element of $\mathcal{V}(\Omega_{\varepsilon})$ satisfying \eqref{vecje}.
\end{proof}

We give the following definition for the weak solutions of problem \eqref{nsstokespd1} (equally, of problem \eqref{nsstokespd}):
\begin{definition}\label{weaksolution}
A vector field $u \in \mathcal{V}_{\sigma}(\Omega_{\varepsilon})$ is called a \textbf{weak solution} of \eqref{nsstokespd1} if
$$
		\int_{\Omega_{\varepsilon}} \nabla u \cdot \nabla \varphi + \int_{\Omega_{\varepsilon}} \left[ \nabla u -  (\nabla u)^{\top} \right] u \cdot \varphi + p^{+} \left( \int_{\Gamma_{O}} \varphi \cdot \widehat{k} \right) - p^{-} \left( \int_{\Gamma_{I}} \varphi \cdot \widehat{k} \right) = \int_{\Omega_{\varepsilon}} f \cdot \varphi \qquad \forall \varphi \in \mathcal{V}_{\sigma}(\Omega_{\varepsilon}) \, .
$$
\end{definition}
\noindent
For the sake of clarity, let us explain how the weak formulation of Definition \ref{weaksolution} arises. Suppose that $u=(u_1,u_2,u_3) \in \mathcal{C}^{2}(\overline{\Omega_{\varepsilon}})$ and $\Phi \in \mathcal{C}^{1}(\overline{\Omega_{\varepsilon}})$ solve \eqref{nsstokespd1} in the classical sense. Given any vector field $\varphi=(\varphi_1,\varphi_2,\varphi_3) \in H^{1}(\Omega_{\varepsilon})$, we integrate by parts (separately in $\Omega_{\pm}$) in the following way:
\begin{equation} \label{testing1}
\begin{aligned}
& \int_{\Omega_{\varepsilon}} ( -\Delta u + \nabla \Phi) \cdot \varphi = \int_{\Omega_{-}} ( -\Delta u^{-} + \nabla \Phi^{-}) \cdot \varphi^{-} + \int_{\Omega_{+}} ( -\Delta u^{+} + \nabla \Phi^{+})  \cdot \varphi^{+} \\[6pt]
& \hspace{-5mm} = \int_{\Omega_{\varepsilon}} \nabla u \cdot \nabla \varphi - \int_{\Omega_{\varepsilon}} \Phi (\nabla \cdot \varphi) + \int_{\partial \Omega_{-}} \left( \Phi^{-}  \nu - \dfrac{\partial u^{-}}{\partial \nu} \right) \cdot \varphi^{-} + \int_{\partial \Omega_{+}} \left( \Phi^{+}  \nu - \dfrac{\partial u^{+}}{\partial \nu} \right) \cdot \varphi^{+} \, .
\end{aligned}
\end{equation}
If, in addition, we assume that $\varphi$ is divergence-free and vanishes on $\Gamma^{\varepsilon}_{W}$, from \eqref{testing1} we obtain
\begin{equation} \label{testing2}
\begin{aligned}
\int_{\Omega_{\varepsilon}} ( -\Delta u + \nabla \Phi) \cdot \varphi & = \int_{\Omega_{\varepsilon}} \nabla u \cdot \nabla \varphi + \int_{\Gamma_{I}} \left( \Phi^{-}  \nu - \dfrac{\partial u^{-}}{\partial \nu} \right) \cdot \varphi^{-} + \int_{K_{\varepsilon}} \left( \Phi^{-}  \nu - \dfrac{\partial u^{-}}{\partial \nu} \right) \cdot \varphi^{-} \\[6pt]
& \hspace{4mm} + \int_{\Gamma_{O}} \left( \Phi^{+}  \nu - \dfrac{\partial u^{+}}{\partial \nu} \right) \cdot \varphi^{+} + \int_{K_{\varepsilon}} \left( \Phi^{+}  \nu - \dfrac{\partial u^{+}}{\partial \nu} \right) \cdot \varphi^{+} \\[6pt]
& = \int_{\Omega_{\varepsilon}} \nabla u \cdot \nabla \varphi + \int_{\Gamma_{I}} \left( \Phi^{-}  \nu - \dfrac{\partial u^{-}}{\partial \nu} \right) \cdot \varphi^{-} + \int_{\Gamma_{O}} \left( \Phi^{+}  \nu - \dfrac{\partial u^{+}}{\partial \nu} \right) \cdot \varphi^{+} \, ,
\end{aligned}
\end{equation}
where the integrals over $K_{\varepsilon}$ cancel each other because of the opposite signs of the outward unit normals to $\partial \Omega_{\pm}$ on $\Sigma$ (see Remark \ref{nonlip}). Notice that $\nu = \mp \widehat{k}$ on $\Gamma_{I}$ and $\Gamma_{O}$, respectively, thus $u_{1}=u_{2}=0$ on $\Gamma_{I} \cup \Gamma_{O}$, in view of \eqref{nsstokespd1}$_3$-\eqref{nsstokespd1}$_4$. The regularity and incompressibility condition of $u$ then imply that
$$
\dfrac{\partial u^{-}_{3}}{\partial z} = - \left( \dfrac{\partial u^{-}_{1}}{\partial x} + \dfrac{\partial u^{-}_{2}}{\partial y} \right) = 0 \quad \text{on} \quad \Gamma_{I} \, ; \qquad \dfrac{\partial u^{+}_{3}}{\partial z} = - \left( \dfrac{\partial u^{+}_{1}}{\partial x} + \dfrac{\partial u^{+}_{2}}{\partial y} \right) = 0 \quad \text{on} \quad \Gamma_{O} \, .
$$
If we further impose that $\varphi \times \nu = 0$ on $\Gamma_{I} \cup \Gamma_{O}$ (so that $\varphi_{1}=\varphi_{2}=0$ on $\Gamma_{I} \cup \Gamma_{O}$), we get
\begin{equation} \label{phi1}
	\dfrac{\partial u^{-}}{\partial \nu} \cdot \varphi^{-} = - \dfrac{\partial u^{-}_{3}}{\partial z} \, \varphi^{-}_{3} = 0 \quad \text{on} \quad \Gamma_{I} \, ; \qquad \dfrac{\partial u^{+}}{\partial \nu} \cdot \varphi^{+} = \dfrac{\partial u^{+}_{3}}{\partial z} \, \varphi^{+}_{3} = 0 \quad \text{on} \quad \Gamma_{O} \, .
\end{equation}
Inserting \eqref{phi1} into \eqref{testing2}, together with the boundary conditions verified by the Bernoulli pressure in \eqref{nsstokespd1}$_3$-\eqref{nsstokespd1}$_4$, yields the identity
\begin{equation} \label{testing3}
\int_{\Omega_{\varepsilon}} ( -\Delta u + \nabla \Phi) \cdot \varphi = \int_{\Omega_{\varepsilon}} \nabla u \cdot \nabla \varphi + p^{+} \left( \int_{\Gamma_{O}} \varphi \cdot \widehat{k} \right) - p^{-} \left( \int_{\Gamma_{I}} \varphi \cdot \widehat{k} \right) \, .
\end{equation}
Finally, multiplying the equation of conservation of momentum \eqref{nsstokespd1}$_1$ by a vector field $\varphi \in \mathcal{V}_{\sigma}(\Omega_{\varepsilon})$ and integrating by parts in $\Omega_{\pm}$ (enforcing  \eqref{testing3}), provides the weak formulation of Definition \ref{weaksolution}. Conversely, we refer to \cite[Section 2]{korobkov2020solvability} for an explanation of the fact that the boundary conditions involving the Bernoulli pressure in \eqref{nsstokespd1} are implicitly contained in the variational formulation of Definition \ref{weaksolution}.

The main result of this section guarantees the unrestricted solvability (concerning the size of the data) of problem \eqref{nsstokespd1}.
\begin{theorem} \label{epslevel}
	Let $\Omega_{\varepsilon}$ be as in \eqref{perfordomain}. For any $p^{\pm} \in \mathbb{R}$ and $f \in L^{2}(\Omega)$, there exists at least one weak solution $u_{\varepsilon} \in \mathcal{V}_{\sigma}(\Omega_{\varepsilon})$ of the prescribed pressure drop problem \eqref{nsstokespd1} and a unique associated pressure $\Phi_{\varepsilon} \in L^{2}(\Omega_{\varepsilon})$ such that the pair $(u_{\varepsilon},\Phi_{\varepsilon})$ satisfies the following identity:
	\begin{equation} \label{pressurewf}
	\begin{aligned}	
	& \int_{\Omega_{\varepsilon}} \nabla u_{\varepsilon} \cdot \nabla \varphi + \int_{\Omega_{\varepsilon}} \left[ \nabla u_{\varepsilon} -  (\nabla u_{\varepsilon})^{\top} \right] u_{\varepsilon} \cdot \varphi + p^{+} \left( \int_{\Gamma_{O}} \varphi \cdot \widehat{k} \right) \\[6pt]
	& - p^{-} \left( \int_{\Gamma_{I}} \varphi \cdot \widehat{k} \right) - \int_{\Omega_{\varepsilon}} \Phi_{\varepsilon} (\nabla \cdot \varphi) = \int_{\Omega_{\varepsilon}} f \cdot \varphi \qquad \forall \varphi \in \mathcal{V}(\Omega_{\varepsilon}) \, .
	\end{aligned}
	\end{equation}
\end{theorem}
\noindent
\begin{proof}
	Since $\varepsilon \in I_{*}$ is a fixed parameter in this section, the corresponding subscript will be omitted sometimes for the sake of simplicity. Accordingly, throughout this proof, $C > 0$ will denote a generic constant that depends on $\Omega_{\varepsilon}$, but that may change from line to line.
	\par
Given any $\varphi \in \mathcal{V}_{\sigma}(\Omega_{\varepsilon})$, by arguing as in \eqref{refer00} one can easily prove that
\begin{equation} \label{sameflux}
\int_{\Gamma_{I}} \varphi \cdot \widehat{k} = \int_{\Gamma_{O}} \varphi \cdot \widehat{k} \, .
\end{equation}
In view of Definition \ref{weaksolution}, we then search for a vector field $u \in \mathcal{V}_{\sigma}(\Omega_{\varepsilon})$ such that
	\begin{equation} \label{oseenquasipd}
		\int_{\Omega_{\varepsilon}} \nabla u \cdot \nabla \varphi + \int_{\Omega_{\varepsilon}} \mathcal{E}(u) u \cdot \varphi + (p^{+}-p^{-}) \int_{\Gamma_{O}} \varphi \cdot \widehat{k} = \int_{\Omega_{\varepsilon}} f \cdot \varphi \qquad \forall \varphi \in \mathcal{V}_{\sigma}(\Omega_{\varepsilon}) \, .
	\end{equation}
In \eqref{oseenquasipd} we have denoted by $\mathcal{E}(w) \doteq \nabla w - (\nabla w)^{\top}$ the skew-symmetric gradient of any vector field $w \in H^{1}(\Omega_{\varepsilon})$. For a given and fixed $u \in \mathcal{V}_{\sigma}(\Omega_{\varepsilon})$, the applications
	$$
	\varphi \in \mathcal{V}_{\sigma}(\Omega_{\varepsilon}) \longmapsto \int_{\Omega_{\varepsilon}} \mathcal{E}(u)u \cdot \varphi \qquad \text{and} \qquad
	\varphi \in \mathcal{V}_{\sigma}(\Omega_{\varepsilon}) \longmapsto  \int_{\Omega_{\varepsilon}} f \cdot \varphi - (p^{+}-p^{-}) \int_{\Gamma_{O}} \varphi \cdot \widehat{k}
	$$
	clearly define linear continuous functions on $\mathcal{V}_{\sigma}(\Omega_{\varepsilon})$. Then, in view of the Riesz Representation Theorem, the identity \eqref{oseenquasipd} may be written as
	$$
	[u + \mathcal{P}(u) - \mathcal{F}, \varphi]_{\mathcal{V}_{\sigma}(\Omega_{\varepsilon})} = 0 \qquad \forall \varphi \in \mathcal{V}_{\sigma}(\Omega_{\varepsilon}) \, ,
	$$
	see \eqref{dirscal}, for some (unique) elements $\mathcal{P}(u), \mathcal{F} \in \mathcal{V}_{\sigma}(\Omega_{\varepsilon})$ such that
	$$
	[\mathcal{P}(u), \varphi]_{\mathcal{V}(\Omega_{\varepsilon})} = \int_{\Omega_{\varepsilon}} \mathcal{E}(u)u \cdot \varphi \quad \text{and} \quad
	[\mathcal{F}, \varphi]_{\mathcal{V}(\Omega_{\varepsilon})} = \int_{\Omega_{\varepsilon}} f \cdot \varphi - (p^{+}-p^{-}) \int_{\Gamma_{O}} \varphi \cdot \widehat{k} \qquad \forall \varphi \in \mathcal{V}_{\sigma}(\Omega_{\varepsilon}) \, .
	$$
	We have so defined a linear operator $\mathcal{P} : \mathcal{V}_{\sigma}(\Omega_{\varepsilon}) \longrightarrow \mathcal{V}_{\sigma}(\Omega_{\varepsilon})$ and we are led to find a solution $u \in \mathcal{V}_{\sigma}(\Omega_{\varepsilon})$ of the following nonlinear operator equation:
	\begin{equation}\label{oseenquasipd2}
		u + \mathcal{P}(u) - \mathcal{F} = 0 \quad \text{in} \ \ \mathcal{V}_{\sigma}(\Omega_{\varepsilon}) \, .
	\end{equation}
	Exactly as in \cite[Chapter 5, Theorem 1]{ladyzhenskaya1969mathematical} one can show that the operator $\mathcal{P}$ is compact. Therefore, as a consequence of the Leray-Schauder Principle \cite[Chapter 6]{zeidler2013nonlinear}, in order to prove that \eqref{oseenquasipd2} possesses at least one solution, it suffices to guarantee that any $v^{\lambda} \in \mathcal{V}_{\sigma}(\Omega_{\varepsilon})$ such that
	\begin{equation}\label{oseenquasipd3}
		v^{\lambda} + \lambda (\mathcal{P}(v^{\lambda}) - \mathcal{F}) = 0 \quad \text{in} \ \ \mathcal{V}_{\sigma}(\Omega_{\varepsilon}) \, ,
	\end{equation}
	is uniformly bounded with respect to $\lambda \in [0, 1]$. Given $\lambda \in [0, 1]$ and $v^{\lambda} \in \mathcal{V}_{\sigma}(\Omega_{\varepsilon}) \setminus \{0\}$ such that \eqref{oseenquasipd3} holds, we clearly have
	\begin{equation} \label{vlambdapd0}
		\int_{\Omega_{\varepsilon}} \nabla v^{\lambda} \cdot \nabla \varphi + \lambda \int_{\Omega_{\varepsilon}} \mathcal{E}(v^{\lambda})v^{\lambda} \cdot \varphi + \lambda (p^{+}-p^{-}) \int_{\Gamma_{O}} \varphi \cdot \widehat{k} = \lambda \int_{\Omega_{\varepsilon}} f \cdot \varphi \qquad \forall \varphi \in \mathcal{V}_{\sigma}(\Omega_{\varepsilon}) \, .
	\end{equation}
	By putting $\varphi = v^{\lambda}$ in \eqref{vlambdapd0} we obtain
	\begin{equation} \label{vlambdapd1}
		\| \nabla v^{\lambda} \|^{2}_{L^{2}(\Omega_{\varepsilon})} = \lambda  \int_{\Omega_{\varepsilon}} f \cdot v^{\lambda} - \lambda (p^{+}-p^{-}) \int_{\Gamma_{O}} v^{\lambda} \cdot \widehat{k} \, .
	\end{equation}
	Applying the H\"older, Poincaré and trace inequalities, in $\Omega_{\varepsilon}$ (see Lemma \ref{unisob}), on the right-hand side of \eqref{vlambdapd1} gives us the bound
	$$
    \| \nabla v^{\lambda} \|_{L^{2}(\Omega_{\varepsilon})} \leq C \left( \| f \|_{L^{2}(\Omega_{\varepsilon})} + | p^{+}-p^{-} |  \right) \qquad \forall \lambda \in [0,1] \, ,
	$$
	which shows that problem \eqref{nsstokespd1} admits, at least, one weak solution $u_{\varepsilon} \in \mathcal{V}_{\sigma}(\Omega_{\varepsilon})$. 

Denote by $\mathcal{V}^{-1}(\Omega_{\varepsilon})$ the dual space of $\mathcal{V}(\Omega_{\varepsilon})$ and by $\langle\cdot,\cdot\rangle_{\mathcal{V}(\Omega_{\varepsilon})}$ the
	duality product
	between $\mathcal{V}^{-1}(\Omega_{\varepsilon})$ and $\mathcal{V}(\Omega_{\varepsilon})$. The gradient
	$-\text{grad}_{\Omega_{\varepsilon}} : L^{2}(\Omega_{\varepsilon}) \longrightarrow \mathcal{V}^{-1}(\Omega_{\varepsilon})$ of a scalar function $q \in L^{2}(\Omega_{\varepsilon})$ is defined by
	$$
	\langle - \text{grad}_{\Omega_{\varepsilon}}(q) , \varphi \rangle_{\mathcal{V}(\Omega_{\varepsilon})} \doteq \int_{\Omega_{\varepsilon}} q (\nabla \cdot \varphi) \qquad \forall \varphi \in \mathcal{V}(\Omega_{\varepsilon})\, ,
	$$
	thus being the adjoint of the (strong) divergence operator $\text{div}_{\Omega_{\varepsilon}} : \mathcal{V}(\Omega_{\varepsilon}) \longrightarrow L^{2}(\Omega_{\varepsilon})$.
	Therefore, the Closed Range Theorem of Banach \cite[Chapter VII]{yosida2012functional} can be applied to deduce that
	$$
	\text{Range}(-\text{grad}_{\Omega_{\varepsilon}}) = \text{Ker}(\text{div}_{\Omega_{\varepsilon}})^{\perp} = \mathcal{V}_{\sigma}(\Omega_{\varepsilon})^{\perp} \doteq \{ \mathcal{F} \in
	\mathcal{V}^{-1}(\Omega_{\varepsilon}) \ | \ \langle \mathcal{F} , \varphi \rangle_{\mathcal{V}(\Omega_{\varepsilon})} = 0 \quad \forall \varphi \in
	\mathcal{V}_{\sigma}(\Omega_{\varepsilon}) \, \} \, .
	$$
	We then define $\mathcal{Q}_{\varepsilon} \in \mathcal{V}^{-1}(\Omega_{\varepsilon})$ by
	$$
	\begin{aligned}
		\mathcal{Q}_{\varepsilon}(\varphi) \doteq \int_{\Omega_{\varepsilon}} \nabla u_{\varepsilon} \cdot \nabla \varphi + \int_{\Omega_{\varepsilon}} \mathcal{E}(u_{\varepsilon}) u_{\varepsilon} \cdot \varphi + p^{+} \left( \int_{\Gamma_{O}} \varphi \cdot \widehat{k} \right) - p^{-} \left( \int_{\Gamma_{I}} \varphi \cdot \widehat{k} \right) - \int_{\Omega_{\varepsilon}} f \cdot \varphi \quad \forall \varphi \in \mathcal{V}(\Omega_{\varepsilon}) \, .
	\end{aligned}	
	$$
In view of \eqref{sameflux}-\eqref{oseenquasipd} we deduce that $\mathcal{Q}_{\varepsilon} \in \text{Range}(-\text{grad}_{\Omega_{\varepsilon}})$, thus guaranteeing the existence of a scalar pressure $\Phi_{\varepsilon} \in L^{2}(\Omega_{\varepsilon})$ such that
\begin{equation} \label{presarise}
\mathcal{Q}_{\varepsilon}(\varphi)  = \int_{\Omega_{\varepsilon}} \Phi_{\varepsilon} (\nabla \cdot \varphi) \qquad \forall \varphi \in \mathcal{V}(\Omega_{\varepsilon}) \, .
\end{equation}
Suppose there exists another scalar function $q_{\varepsilon} \in L^{2}(\Omega_{\varepsilon})$ verifying \eqref{presarise}, thereby implying
\begin{equation} \label{bog1}
	\int_{\Omega_{\varepsilon}} (\Phi_{\varepsilon}-q_{\varepsilon}) (\nabla \cdot \varphi) = 0 \qquad \forall \varphi \in \mathcal{V}(\Omega_{\varepsilon}) \, .
\end{equation}
Lemma \ref{bogtype} provides the existence of a vector field $Y_{\varepsilon} \in \mathcal{V}(\Omega_{\varepsilon})$ such that $\nabla \cdot Y_{\varepsilon} = \Phi_{\varepsilon}-q_{\varepsilon}$ in $\Omega_{\varepsilon}$. Taking $\varphi = Y_{\varepsilon}$ in \eqref{bog1} ensures that $\Phi_{\varepsilon}=q_{\varepsilon}$ almost everywhere in $\Omega_{\varepsilon}$. This concludes the proof
\end{proof}

\begin{remark} \label{additiveconstant}
In virtue of the Divergence Theorem we observe that
$$
\int_{\Omega_{\varepsilon}} \nabla \cdot \varphi = \int_{\Gamma_{O}} \varphi \cdot \widehat{k} - \int_{\Gamma_{I}} \varphi \cdot \widehat{k} \qquad \forall \varphi \in \mathcal{V}(\Omega_{\varepsilon}) \, .
$$
Thus, the (unique) Bernoulli pressure associated to any weak solution of problem \eqref{nsstokespd1} can be modified, up to an additive constant, at the expense of altering the prescribed pressure drop. In particular, given a weak solution $u_{\varepsilon} \in \mathcal{V}_{\sigma}(\Omega_{\varepsilon})$ of \eqref{nsstokespd1}, its associated pressure $\Phi_{\varepsilon} \in L^{2}(\Omega_{\varepsilon})$ can be recast in such a way that $\Phi_{\varepsilon} \in L_{0}^{2}(\Omega_{\varepsilon})$, at the cost of introducing an $\varepsilon$-dependent term in the pressure drop, namely, the mean value of $\Phi_{\varepsilon}$ over $\Omega_{\varepsilon}$. That is the reason why we do not project the Bernoulli pressure onto $L_{0}^{2}(\Omega_{\varepsilon})$.
\end{remark}

\begin{remark} \label{wfpressure}
In view of Theorem \ref{epslevel}, in the sequel, any pair $(u,\Phi) \in \mathcal{V}_{\sigma}(\Omega_{\varepsilon}) \times L^{2}(\Omega_{\varepsilon})$ will be called a \textbf{weak solution} of problem \eqref{nsstokespd1} if it satisfies the integral identity \eqref{pressurewf} for every $\varphi \in \mathcal{V}(\Omega_{\varepsilon})$.
\end{remark}

\section{Boundary-value problem at the $\varepsilon$-level: uniform bounds} \label{epslevelsec2}

Let $\varepsilon \in I_{*}$ be a fixed parameter. Denoting $\mathcal{L}_{\pm} \doteq \mathcal{L} \cap \partial \Omega_{\pm}$, we further introduce the following two spaces of vector fields:
\begin{equation} \label{masespacios}
\begin{aligned}
& \mathcal{V}(\Omega_{-}) \doteq \left\lbrace v \in H^{1}(\Omega_{-}) \ | \ v \times \nu = 0 \ \ \mbox{on} \ \ \Gamma_{I} \, ; \qquad v = 0 \ \ \mbox{on} \ \ \mathcal{L}_{-} \cup \Sigma \, \right\rbrace \, , \\[6pt]
& \mathcal{V}(\Omega_{+}) \doteq \left\lbrace v \in H^{1}(\Omega_{+}) \ | \ v \times \nu = 0 \ \ \mbox{on} \ \ \Gamma_{O} \, ; \qquad v = 0 \ \ \mbox{on} \ \ \mathcal{L}_{+} \cup \Sigma \, \right\rbrace \, .
\end{aligned}
\end{equation}
\noindent
The main result of this section provides uniform bounds (with respect to $\varepsilon \in I_{*}$) for the solutions of problem \eqref{nsstokespd1}.

\begin{theorem} \label{epslevelbounds}
	Let $\Omega_{\varepsilon}$ be as in \eqref{perfordomain}. For any $p^{\pm} \in \mathbb{R}$ and $f \in L^{2}(\Omega)$, let $(u_{\varepsilon},\Phi_{\varepsilon}) \in \mathcal{V}_{\sigma}(\Omega_{\varepsilon}) \times L^{2}(\Omega_{\varepsilon})$ be a weak solution of problem \eqref{nsstokespd1}. Then, the uniform bound
	\begin{equation} \label{uboundpf}
		\sup_{\varepsilon \in I_{*}} \left( \| \nabla u_{\varepsilon} \|_{L^{2}(\Omega_{\varepsilon})} + \| \Phi_{\varepsilon} \|_{L^{2}(\Omega_{\varepsilon})} \right) \leq C_{*} \, ,
	\end{equation}
	holds for some constant $C_{*} > 0$ that depends on $\Omega$, $p^{\pm}$, $f$ and $\{ \delta_{0}, \delta_{1} \}$, but is independent of $\varepsilon \in I_{*}$.
\end{theorem}
\noindent
\begin{proof}
	In what follows, $C > 0$ will always denote a generic constant that depends on $\Omega$ and $\{ \delta_{0}, \delta_{1} \}$ (independently of $\varepsilon \in I_{*}$), but that may change from line to line.
	\par
We start by noticing that, in view of \eqref{refer00}, we have
$$
\int_{\Gamma_{I}} u_{\varepsilon} \cdot \widehat{k} = \int_{\Gamma_{O}} u_{\varepsilon} \cdot \widehat{k} = \int_{\Sigma} u_{\varepsilon} \cdot \widehat{k} \, .
$$	
Then, by taking $\varphi = u_{\varepsilon}$ as a test function in the weak formulation of Definition \ref{weaksolution} we get
\begin{equation} \label{ueps1}
	\| \nabla u_{\varepsilon}  \|^{2}_{L^{2}(\Omega_{\varepsilon})} = \int_{\Omega_{\varepsilon}} f \cdot u_{\varepsilon} - (p^{+}-p^{-}) \int_{\Sigma} u_{\varepsilon}  \cdot \widehat{k} \, .
\end{equation}
Applying the H\"older and Poincaré inequalities (see Lemma \ref{unisob}), together with \eqref{conca0}, on the right-hand side of \eqref{ueps1}, gives us
$$
\begin{aligned}
\| \nabla u_{\varepsilon}  \|^{2}_{L^{2}(\Omega_{\varepsilon})} & \leq \| f  \|_{L^{2}(\Omega_{\varepsilon})} \| u_{\varepsilon}  \|_{L^{2}(\Omega_{\varepsilon})} + \sqrt{\pi} \, R \, | p^{+}-p^{-} | \,  \| u_{\varepsilon}  \|_{L^{2}(\Sigma)} \\[6pt]
& \leq C \left( \| f  \|_{L^{2}(\Omega)}  + \sqrt{r_{\varepsilon}} \, | p^{+}-p^{-} | \right) \| \nabla u_{\varepsilon} \|_{L^{2}(\Omega_{\varepsilon})} \, ,
\end{aligned}
$$
so that
\begin{equation} \label{ueps2}
\| \nabla u_{\varepsilon}  \|_{L^{2}(\Omega_{\varepsilon})} \leq C \left( \| f  \|_{L^{2}(\Omega)}  + | p^{+}-p^{-} | \right) \qquad \forall \varepsilon \in I_{*} \, .
\end{equation}

Now, given any vector field $\varphi \in \mathcal{V}(\Omega_{\pm})$, by trivial extension to the complement $\Omega_{\mp}$ we have that $\varphi \in \mathcal{V}(\Omega_{\varepsilon})$. Once taken as a test function in the weak formulation of \eqref{nsstokespd1} (see Remark \ref{wfpressure}), this gives 
\begin{equation} \label{split1}
\int_{\Omega_{-}} \nabla u^{-}_{\varepsilon} \cdot \nabla \varphi + \int_{\Omega_{-}} \mathcal{E}(u^{-}_{\varepsilon}) u^{-}_{\varepsilon} \cdot \varphi - p^{-} \left( \int_{\Gamma_{I}} \varphi \cdot \widehat{k} \right) - \int_{\Omega_{-}} \Phi^{-}_{\varepsilon} (\nabla \cdot \varphi) = \int_{\Omega_{-}} f^{-} \cdot \varphi \quad \forall \varphi \in \mathcal{V}(\Omega_{-}) \, ,
\end{equation}
and also
\begin{equation} \label{split2}
\int_{\Omega_{+}} \nabla u^{+}_{\varepsilon} \cdot \nabla \varphi + \int_{\Omega_{+}} \mathcal{E}(u^{+}_{\varepsilon}) u^{+}_{\varepsilon} \cdot \varphi + p^{+} \left( \int_{\Gamma_{O}} \varphi \cdot \widehat{k} \right) - \int_{\Omega_{+}} \Phi^{+}_{\varepsilon} (\nabla \cdot \varphi) = \int_{\Omega_{+}} f^{+} \cdot \varphi \quad \forall \varphi \in \mathcal{V}(\Omega_{+}) \, .
\end{equation}
On the other hand, the Divergence Theorem trivially implies that
\begin{equation} \label{split3}
\int_{\Omega_{-}} \nabla \cdot \varphi = - \int_{\Gamma_{I}} \varphi \cdot \widehat{k} \qquad \forall \varphi \in \mathcal{V}(\Omega_{-}) \, ; \qquad \int_{\Omega_{+}} \nabla \cdot \varphi = \int_{\Gamma_{O}} \varphi \cdot \widehat{k} \qquad \forall \varphi \in \mathcal{V}(\Omega_{+}) \, .
\end{equation}
Therefore, defining the functions $P^{\pm}_{\varepsilon} \in L_{0}^{2}(\Omega_{\pm})$ by
\begin{equation} \label{projecpres}
P^{\pm}_{\varepsilon} \doteq \Phi^{\pm}_{\varepsilon} - \phi^{\pm}_{\varepsilon} \doteq \Phi^{\pm}_{\varepsilon} - \dfrac{1}{| \Omega_{\pm} |} \int_{\Omega_{\pm}} \Phi^{\pm}_{\varepsilon} \qquad \text{in} \ \ \Omega_{\pm} \, ,
\end{equation}
the identities \eqref{split1}-\eqref{split2} become, owing to \eqref{split3}, 
\begin{equation} \label{split4}
\begin{aligned}	
& \int_{\Omega_{-}} \nabla u^{-}_{\varepsilon} \cdot \nabla \varphi + \int_{\Omega_{-}} \mathcal{E}(u^{-}_{\varepsilon}) u^{-}_{\varepsilon} \cdot \varphi - \int_{\Omega_{-}} P^{-}_{\varepsilon} (\nabla \cdot \varphi) \\[6pt]
& \hspace{-4mm} = \int_{\Omega_{-}} f^{-} \cdot \varphi - \left( \phi^{-}_{\varepsilon} - p^{-} \right) \left( \int_{\Gamma_{I}} \varphi \cdot \widehat{k} \right) \qquad \forall \varphi \in \mathcal{V}(\Omega_{-}) \, ,
\end{aligned}	
\end{equation}
and also
\begin{equation} \label{split5}
\begin{aligned}	
& \int_{\Omega_{+}} \nabla u^{+}_{\varepsilon} \cdot \nabla \varphi + \int_{\Omega_{+}} \mathcal{E}(u^{+}_{\varepsilon}) u^{+}_{\varepsilon} \cdot \varphi  - \int_{\Omega_{+}} P^{+}_{\varepsilon} (\nabla \cdot \varphi) 	\\[6pt]
& \hspace{-4mm}	= \int_{\Omega_{+}} f^{+} \cdot \varphi - \left( p^{+} - \phi^{+}_{\varepsilon} \right) \left( \int_{\Gamma_{O}} \varphi \cdot \widehat{k} \right) \qquad \forall \varphi \in \mathcal{V}(\Omega_{+}) \, .
\end{aligned}		
\end{equation}
In view of \eqref{ueps2}-\eqref{projecpres}, it now suffices to show the uniform bound
\begin{equation} \label{uboundpfproof}
	\sup_{\varepsilon \in I_{*}} \left( \| P^{\pm}_{\varepsilon} \|_{L^{2}(\Omega_{\pm})} + | \phi^{\pm}_{\varepsilon} | \right) \leq C \, .
\end{equation}

Since $P^{\pm}_{\varepsilon} \in L_{0}^{2}(\Omega_{\pm})$ and both domains $\Omega_{\pm}$ have Lipschitz-continuous boundaries, there exist vector fields $J^{\pm}_{\varepsilon} \in H_{0}^{1}(\Omega_{\pm})$ such that
\begin{equation} \label{bogn0}
	\nabla \cdot J^{\pm}_{\varepsilon} = P^{\pm}_{\varepsilon} \quad \text{in} \quad \Omega_{\pm} \qquad \text{and} \qquad  \| \nabla J^{\pm}_{\varepsilon} \|_{L^{2}(\Omega_{\pm})} \leq C \| P^{\pm}_{\varepsilon} \|_{L^{2}(\Omega_{\pm})} \, ,
\end{equation}
see again \cite{bogovskii1979solution}. As $H_{0}^{1}(\Omega_{\pm}) \subset \mathcal{V}(\Omega_{\pm})$, we can take $\varphi = J^{\pm}_{\varepsilon}$ as test functions in \eqref{split4}-\eqref{split5}, thereby obtaining
$$
\int_{\Omega_{\pm}} \nabla u^{\pm}_{\varepsilon} \cdot \nabla J^{\pm}_{\varepsilon} + \int_{\Omega_{\pm}} \mathcal{E}(u^{\pm}_{\varepsilon}) u^{\pm}_{\varepsilon} \cdot J^{\pm}_{\varepsilon} - \| P^{\pm}_{\varepsilon} \|^{2}_{L^{2}(\Omega_{\pm})} = \int_{\Omega_{\pm}} f^{\pm} \cdot J^{\pm}_{\varepsilon} \, .
$$
We then apply the inequalities of H\"older, Poincaré and Sobolev (in $\Omega_{\pm}$), together with \eqref{bogn0}, to estimate
$$
\begin{aligned}
	\| P^{\pm}_{\varepsilon} \|^{2}_{L^{2}(\Omega_{\pm})} & = \int_{\Omega_{\pm}} \nabla u^{\pm}_{\varepsilon} \cdot \nabla J^{\pm}_{\varepsilon} + \int_{\Omega_{\pm}} \mathcal{E}(u^{\pm}_{\varepsilon}) u^{\pm}_{\varepsilon} \cdot J^{\pm}_{\varepsilon} -  \int_{\Omega_{\pm}} f^{\pm} \cdot J^{\pm}_{\varepsilon} \\[6pt]
	& \leq \| \nabla u^{\pm}_{\varepsilon} \|_{L^{2}(\Omega_{\pm})} \| \nabla J^{\pm}_{\varepsilon} \|_{L^{2}(\Omega_{\pm})} + \| \mathcal{E}(u^{\pm}_{\varepsilon}) \|_{L^{2}(\Omega_{\pm})} \| u^{\pm}_{\varepsilon} \|_{L^{4}(\Omega_{\pm})} \| J^{\pm}_{\varepsilon} \|_{L^{4}(\Omega_{\pm})} + \| f \|_{L^{2}(\Omega)} \| J^{\pm}_{\varepsilon} \|_{L^{2}(\Omega_{\pm})} \\[6pt]
	& \leq \| \nabla u^{\pm}_{\varepsilon} \|_{L^{2}(\Omega_{\pm})} \| \nabla J^{\pm}_{\varepsilon} \|_{L^{2}(\Omega_{\pm})} + C \| \nabla u^{\pm}_{\varepsilon} \|^{2}_{L^{2}(\Omega_{\pm})} \| \nabla J^{\pm}_{\varepsilon} \|_{L^{2}(\Omega_{\pm})} + C \| f \|_{L^{2}(\Omega)} \| \nabla J^{\pm}_{\varepsilon} \|_{L^{2}(\Omega_{\pm})} \\[6pt]
	& \leq C \left( 1 + \| \nabla u^{\pm}_{\varepsilon} \|^{2}_{L^{2}(\Omega_{\pm})} + \| f \|_{L^{2}(\Omega)} \right) \| P^{\pm}_{\varepsilon} \|_{L^{2}(\Omega_{\pm})}  \, ,
\end{aligned}
$$
and thus, after further imposing \eqref{ueps2}, we get
\begin{equation} \label{estpresind}
		\| P^{\pm}_{\varepsilon} \|_{L^{2}(\Omega_{\pm})} \leq C \left( 1 + \| f \|^{2}_{L^{2}(\Omega)} + | p^{+}-p^{-} |^{2} \right) \qquad \forall \varepsilon \in I_{*} \, .
\end{equation}
Define the vector field $\varphi_{*} \in \mathcal{C}^{\infty}(\overline{\Omega})$ by
\begin{equation} \label{hagen}
\varphi_{*}(x,y,z) \doteq \dfrac{2}{\pi R^{4}} \left( \dfrac{z}{h} \right)^{2} (R^2 - x^2 - y^2) \widehat{k} \qquad \forall (x,y,z) \in \overline{\Omega} \, ,
\end{equation}
so that clearly $\varphi_{*}^{\pm} \in \mathcal{V}(\Omega_{\pm})$ with
$$
\int_{\Gamma_{I}} \varphi_{*}^{-} \cdot \widehat{k} = \int_{\Gamma_{O}} \varphi_{*}^{+} \cdot \widehat{k} = 1 \, .
$$
Taking $\varphi = \varphi_{*}^{\pm}$ as test functions in \eqref{split4}-\eqref{split5} we obtain, respectively, 
\begin{equation} \label{split6}
\phi^{-}_{\varepsilon} = p^{-} + \int_{\Omega_{-}} f^{-} \cdot \varphi_{*}^{-} - \int_{\Omega_{-}} \nabla u^{-}_{\varepsilon} \cdot \nabla \varphi_{*}^{-} - \int_{\Omega_{-}} \mathcal{E}(u^{-}_{\varepsilon}) u^{-}_{\varepsilon} \cdot \varphi_{*}^{-} + \int_{\Omega_{-}} P^{-}_{\varepsilon} (\nabla \cdot \varphi_{*}^{-}) \, ,
\end{equation}
and also
\begin{equation} \label{split7}
\phi^{+}_{\varepsilon} = p^{+} - \int_{\Omega_{+}} f^{+} \cdot \varphi_{*}^{+} + \int_{\Omega_{+}} \nabla u^{+}_{\varepsilon} \cdot \nabla \varphi_{*}^{+} + \int_{\Omega_{+}} \mathcal{E}(u^{+}_{\varepsilon}) u^{+}_{\varepsilon} \cdot \varphi_{*}^{+} - \int_{\Omega_{+}} P^{+}_{\varepsilon} (\nabla \cdot \varphi_{*}^{+})  \, .
\end{equation}
Then, applying H\"older's inequality and \eqref{sob0} to both identities \eqref{split6}-\eqref{split7} gives us the bounds
$$
| \phi^{\pm}_{\varepsilon} | \leq C \left( | p^{\pm} | + \| f \|_{L^{2}(\Omega)} + \| \nabla u^{\pm}_{\varepsilon} \|^{2}_{L^{2}(\Omega_{\pm})}  + \| P^{\pm}_{\varepsilon} \|_{L^{2}(\Omega_{\pm})} \right) \, ,
$$
and again from \eqref{ueps2}-\eqref{estpresind} we conclude
\begin{equation} \label{split8}
| \phi^{\pm}_{\varepsilon} | \leq C \left( 1 + | p^{\pm} |^{2} + \| f \|^{2}_{L^{2}(\Omega)} \right) \qquad \forall \varepsilon \in I_{*} \, .
\end{equation}
Putting together \eqref{estpresind}-\eqref{split8} provides the uniform bound \eqref{uboundpfproof}, thus finishing the proof.
\end{proof}

\begin{remark} \label{bogfails}
Let $\varepsilon \in I_{*}$ and $(u_{\varepsilon},\Phi_{\varepsilon}) \in \mathcal{V}_{\sigma}(\Omega_{\varepsilon}) \times L^{2}(\Omega_{\varepsilon})$ be a weak solution of problem \eqref{nsstokespd1}. From Lemma \ref{bogtype} we deduce the existence of a vector field $Y_{\varepsilon} \in \mathcal{V}(\Omega_{\varepsilon})$ such that
\begin{equation} \label{vecjepre}
	\nabla \cdot Y_{\varepsilon} = \Phi_{\varepsilon} \quad \text{in} \quad \Omega_{\varepsilon} \qquad \text{and} \qquad  \| \nabla Y_{\varepsilon} \|_{L^{2}(\Omega_{\varepsilon})} \leq C_{*} (1 + C_{B}(\Omega_{\varepsilon})) \| \Phi_{\varepsilon} \|_{L^{2}(\Omega_{\varepsilon})} \, ,
\end{equation}
where $C_{B}(\Omega_{\varepsilon}) > 0$ is the Bogovskii constant of $\Omega_{\varepsilon}$ (as in Lemma \ref{bogeps}), and $C_{*} > 0$ is a constant depending on $\Omega$ and $\{ \delta_{0}, \delta_{1} \}$ (independently of $\varepsilon \in I_{*}$). Then, arguing as in the proof of Theorem \eqref{epslevelbounds} (see the estimates prior to \eqref{estpresind}) we can derive the bound
$$
\begin{aligned}
	\| \Phi_{\varepsilon} \|_{L^{2}(\Omega_{\varepsilon})} \leq C_{*} (1 + C_{B}(\Omega_{\varepsilon})) \left( 1 + \| \nabla u_{\varepsilon} \|^{2}_{L^{2}(\Omega_{\varepsilon})} + \| f \|_{L^{2}(\Omega)} \right) \qquad \forall \varepsilon \in I_{*} \, .
\end{aligned}
$$
Nevertheless, in view of \eqref{bogepsineq}-\eqref{ueps2}, this last inequality prevents us from reaching any conclusion concerning the growth of the $L^{2}(\Omega_{\varepsilon})$-norm of the Bernoulli pressure $\Phi_{\varepsilon}$ with respect to $\varepsilon \in I_{*}$. That is the reason why, in the proof of Theorem \eqref{epslevelbounds}, we separated the analysis into each sub-domain $\Omega_{\pm}$.
\end{remark}

Combined with Lemma \ref{concalem}, the uniform bounds obtained in Theorem \eqref{epslevelbounds} allow us to state the following decay rate for the tranverse flux rate and trace over $\Sigma$ of any weak solution of problem \eqref{nsstokespd1}. More precisely:
\begin{corollary} \label{concacor}
	Let $\Omega_{\varepsilon}$ be as in \eqref{perfordomain}. For any $p^{\pm} \in \mathbb{R}$ and $f \in L^{2}(\Omega)$, let $u_{\varepsilon} \in \mathcal{V}_{\sigma}(\Omega_{\varepsilon})$ be a weak solution of problem \eqref{nsstokespd1}, whose tranverse flux rate along $\Omega_{\varepsilon}$ is defined as in \eqref{refer00}. Then, the estimate
\begin{equation} \label{tracezero}
\| u_{\varepsilon}  \|_{L^{2}(\Sigma)} + |F_{\varepsilon}| \leq C_{*} \sqrt{r_{\varepsilon}} \qquad \forall \varepsilon \in I_{*} \, ,
\end{equation}
holds for some constant $C_{*} > 0$ that depends on $\Omega$, $p^{\pm}$, $f$ and $\{ \delta_{0}, \delta_{1} \}$, but is independent of $\varepsilon \in I_{*}$.
\end{corollary}
\noindent
\begin{proof}
	In what follows, $C > 0$ will always denote a generic constant that depends on $\Omega$ and $\{ \delta_{0}, \delta_{1} \}$ (independently of $\varepsilon \in I_{*}$), but that may change from line to line.
	\par
	It suffices to observe that, since $u_{\varepsilon} \in \mathcal{V}_{\sigma}(\Omega_{\varepsilon})$ vanishes on $\Gamma_{\varepsilon}$, as a consequence of H\"older's inequality, the trace inequality \eqref{conca0} and the uniform bound \eqref{uboundpf}, we deduce
	$$
	|F_{\varepsilon}| = \left| \int_{\Sigma} u_{\varepsilon} \cdot \widehat{k} \right| \leq \| u_{\varepsilon} \|_{L^1(\Sigma)} \leq \sqrt{\pi} \, R \, \| u_{\varepsilon} \|_{L^2(\Sigma)} \leq C \sqrt{r_{\varepsilon}} \qquad \forall \varepsilon \in I_{*} \, ,
	$$
	which concludes the proof.
\end{proof}

\section{Asymptotic behavior as $\varepsilon \to 0^{+}$: homogenized equations} \label{energymethod}
By employing the uniform bounds obtained in Theorem \eqref{epslevelbounds}, through a compactness argument (which can be seen as a simplified version of the renowned \textit{energy method} of Tartar \cite[Appendix]{sanchez1980non}) we derive in this section the effective (or \textit{homogenized}) equations satisfied by the solutions of \eqref{nsstokespd1} as $\varepsilon \to 0^{+}$. 

\begin{theorem} \label{effectiveq1}
	Let $(\Omega_{\varepsilon})_{\varepsilon \in I_{*}}$ be the family of domains \eqref{perfordomain}. Given any $p^{\pm} \in \mathbb{R}$ and $f \in L^{2}(\Omega)$, let $(u_{\varepsilon},\Phi_{\varepsilon}) \in \mathcal{V}_{\sigma}(\Omega_{\varepsilon}) \times L^{2}(\Omega_{\varepsilon})$ be a weak solution of problem \eqref{nsstokespd1}. Then, up to the extraction of a subsequence, the sequences of restrictions $\{(u^{\pm}_{\varepsilon}, \Phi^{\pm}_{\varepsilon}) \}_{\varepsilon \in I_{*}} \subset H^{1}(\Omega_{\pm}) \times L^{2}(\Omega_{\pm})$ converge strongly to weak solutions $(u^{\pm}, \Phi^{\pm}) \in H^{1}(\Omega_{\pm}) \times L^{2}(\Omega_{\pm})$ of the following boundary-value problems, in $\Omega_{\pm}$, as $\varepsilon \to 0^{+}$:
	\begin{equation}\label{nsstokespdfinal-}
		\left\{
		\begin{aligned}
			& -\Delta u^{-} + (u^{-}\cdot\nabla)u^{-} - (\nabla u^{-})^{\top} u^{-} + \nabla \Phi^{-} = f^{-} \, ,\ \quad  \nabla\cdot u^{-} = 0 \ \ \mbox{ in } \ \ \Omega_{-} \, , \\[3pt]
			& u^{-} \times \nu = 0 \, , \quad \Phi^{-} = p^{-} \ \ \mbox{ on } \ \ \Gamma_{I} \, ,\\[3pt]
			& u^{-}=0 \ \ \mbox{ on } \ \ \mathcal{L}_{-} \cup \Sigma \, ,
		\end{aligned}
		\right.
	\end{equation}
	and
	\begin{equation}\label{nsstokespdfinal+}
		\left\{
	\begin{aligned}
		& -\Delta u^{+} + (u^{+}\cdot\nabla)u^{+} - (\nabla u^{+})^{\top} u^{+} + \nabla \Phi^{+} = f^{+} \, ,\ \quad  \nabla\cdot u^{+} = 0 \ \ \mbox{ in } \ \ \Omega_{+} \, , \\[3pt]
		& u^{+} \times \nu = 0 \, , \quad \Phi^{+} = p^{+} \ \ \mbox{ on } \ \ \Gamma_{O} \, ,\\[3pt]
		& u^{+}=0 \ \ \mbox{ on } \ \ \mathcal{L}_{+} \cup \Sigma \, .
	\end{aligned}
	\right.
	\end{equation}
Furthermore, $(u^{\pm}, \Phi^{\pm}) \in H^{2}(\Omega_{\pm}) \times H^{1}(\Omega_{\pm})$ and they satisfy in strong form the systems \eqref{nsstokespdfinal-}-\eqref{nsstokespdfinal+}, respectively.	
\end{theorem}
\noindent
\begin{proof}
	In what follows, $C > 0$ will always denote a generic constant that depends on $\Omega$ and $\{ \delta_{0}, \delta_{1} \}$ (independently of $\varepsilon \in I_{*}$), but that may change from line to line.
\par
Let us introduce the following (closed) spaces of solenoidal vector fields: 
$$
\begin{aligned}
& \mathcal{Y}_{\sigma}(\Omega_{-}) \doteq \left\lbrace v \in H^{1}(\Omega_{-}) \ | \ \nabla \cdot v=0 \ \ \mbox{in} \ \ \Omega_{-} \, ; \qquad v \times \nu = 0 \ \ \mbox{on} \ \ \Gamma_{I} \, ; \qquad v = 0 \ \ \mbox{on} \ \ \mathcal{L}_{-} \, \right\rbrace \, , \\[6pt]
& \mathcal{Y}_{\sigma}(\Omega_{+}) \doteq \left\lbrace v \in H^{1}(\Omega_{+}) \ | \ \nabla \cdot v=0 \ \ \mbox{in} \ \ \Omega_{+} \, ; \qquad v \times \nu = 0 \ \ \mbox{on} \ \ \Gamma_{O} \, ; \qquad v = 0 \ \ \mbox{on} \ \ \mathcal{L}_{+} \, \right\rbrace \, .
\end{aligned}
$$
The bounds in \eqref{uboundpf} entail that the sequences of restrictions $\{(u^{\pm}_{\varepsilon}, \Phi^{\pm}_{\varepsilon}) \}_{\varepsilon \in I_{*}} \subset \mathcal{Y}_{\sigma}(\Omega_{\pm}) \times L^{2}(\Omega_{\pm})$ are uniformly bounded, so there exist $(u^{\pm}, \Phi^{\pm}) \in \mathcal{Y}_{\sigma}(\Omega_{\pm}) \times L^{2}(\Omega_{\pm})$ such that the following convergences hold as $\varepsilon \to 0^{+}$:
\begin{equation} \label{convergences0}
 u^{\pm}_{\varepsilon} \rightharpoonup u^{\pm} \ \ \text{weakly in} \ \mathcal{Y}_{\sigma}(\Omega_{\pm}) \, ; \quad u^{\pm}_{\varepsilon} \to u^{\pm} \ \ \text{strongly in} \ L^{4}(\Omega_{\pm}) \, ; \quad \Phi^{\pm}_{\varepsilon} \rightharpoonup \Phi^{\pm} \ \ \text{weakly in} \ L^{2}(\Omega_{\pm}) \, ,
\end{equation}
along subsequences that are not being relabeled, see \cite[Theorem 6.2]{necas2011direct}. Moreover, in \eqref{tracezero} we deduce that, along the subsequences \eqref{convergences0},
\begin{equation} \label{convergences1}
\lim\limits_{\varepsilon \to 0^{+}} \| u^{\pm}_{\varepsilon}  \|_{L^{2}(\Sigma)} = 0 \, ,
\end{equation}
so that, in fact, $u^{\pm} \in \mathcal{Y}_{\sigma}(\Omega_{\pm}) \cap \mathcal{V}(\Omega_{\pm})$. Now, given any test function $\varphi \in \mathcal{V}(\Omega_{\pm})$, the identities \eqref{split1}-\eqref{split2} are valid for every $\varepsilon \in I_{*}$, along the subsequences \eqref{convergences0}. With the help of the convergences in \eqref{convergences0} we can easily prove that
\begin{equation} \label{convergences2}
	\begin{aligned}
		& \lim_{\varepsilon \to 0^{+}} \int_{\Omega_{\pm}} \nabla u^{\pm}_{\varepsilon} \cdot \nabla \varphi = \int_{\Omega_{\pm}} \nabla u^{\pm} \cdot \nabla \varphi \, , \qquad  \lim_{\varepsilon \to 0^{+}} \int_{\Omega_{\pm}} \mathcal{E}(u^{\pm}_{\varepsilon})u^{\pm}_{\varepsilon} \cdot \varphi = \int_{\Omega_{\pm}} \mathcal{E}(u^{\pm})u^{\pm}\cdot \varphi \, , \\[6pt]
		& \hspace{3.2cm} \lim_{\varepsilon \to 0^{+}} \int_{\Omega_{\pm}} \Phi^{\pm}_{\varepsilon} (\nabla \cdot \varphi) = \int_{\Omega_{\pm}} \Phi^{\pm} (\nabla \cdot \varphi) \, ,
	\end{aligned}
\end{equation}
thereby implying
\begin{equation} \label{limit1}
	\int_{\Omega_{-}} \nabla u^{-} \cdot \nabla \varphi + \int_{\Omega_{-}} \mathcal{E}(u^{-}) u^{-} \cdot \varphi - p^{-} \left( \int_{\Gamma_{I}} \varphi \cdot \widehat{k} \right) - \int_{\Omega_{-}} \Phi^{-} (\nabla \cdot \varphi) = \int_{\Omega_{-}} f^{-} \cdot \varphi \quad \forall \varphi \in \mathcal{V}(\Omega_{-}) \, ,
\end{equation}
and also
\begin{equation} \label{limit2}
	\int_{\Omega_{+}} \nabla u^{+} \cdot \nabla \varphi + \int_{\Omega_{+}} \mathcal{E}(u^{+}) u^{+} \cdot \varphi + p^{+} \left( \int_{\Gamma_{O}} \varphi \cdot \widehat{k} \right) - \int_{\Omega_{+}} \Phi^{+} (\nabla \cdot \varphi) = \int_{\Omega_{+}} f^{+} \cdot \varphi \quad \forall \varphi \in \mathcal{V}(\Omega_{+}) \, .
\end{equation}
Notice that, by reasoning as in \eqref{refer00}, we immediately deduce that
\begin{equation} \label{refer01}
	\int_{\Sigma(s)} u^{-} \cdot \widehat{k} = 0 \qquad \forall s \in [-h,0] \, ; \qquad \int_{\Sigma(s)} u^{+} \cdot \widehat{k} = 0 \qquad \forall s \in [0,h] \, .
\end{equation}

Now, consider the domains
$$
	\Omega^{*}_{-} \doteq \left\lbrace (x,y,z) \in \Omega \ \Big| \ -h < z < - \dfrac{h}{2} \, \right\rbrace \qquad \text{and} \qquad \Omega^{*}_{+} \doteq \left\lbrace (x,y,z) \in \Omega \ \Big| \ \dfrac{h}{2} < z < h \, \right\rbrace \, .
$$
Fix any $\varepsilon \in I_{*}$ along the subsequences \eqref{convergences0}. The same extension argument of \cite[Theorem 3.2]{sperone2022} (which employs the usual interior regularity results for the steady-state Navier-Stokes equations, see \cite[Theorem IX.5.1]{galdi2011introduction} or \cite[Section 2]{korobkov2020solvability}) allows us to deduce that $(u^{\pm}_{\varepsilon}, \Phi^{\pm}_{\varepsilon}) \in H^{2}(\Omega^{*}_{\pm}) \times H^{1}(\Omega^{*}_{\pm})$, and so
\begin{equation} \label{gammareg1}
\left\{
\begin{aligned}	
	& -\Delta u^{\pm}_{\varepsilon} + \nabla \Phi^{\pm}_{\varepsilon} = f^{\pm} - \mathcal{E}(u^{\pm}_{\varepsilon}) (u^{\pm}_{\varepsilon}) \quad \text{almost everywhere in} \ \Omega^{*}_{\pm} \, ; \\[6pt]
	& \Phi_{\varepsilon} = p^{-} \ \ \mbox{ almost everywhere on} \ \ \Gamma_{I} \, ; \quad \Phi_{\varepsilon} = p^{+} \ \ \mbox{ almost everywhere on} \ \ \Gamma_{O} \, . 
\end{aligned}
\right.	
\end{equation}
The Sobolev embedding $H^{1}(\Omega_{\pm}) \subset L^{6}(\Omega_{\pm})$ also implies that $f^{\pm} - \mathcal{E}(u^{\pm}_{\varepsilon}) (u^{\pm}_{\varepsilon}) \in L^{3/2}(\Omega_{\pm})$. Moreover, from the inequalities of H\"older and Sobolev (applied in $\Omega_{\pm}$) we derive the estimate
\begin{equation} \label{regsto1}
	\begin{aligned}
		& \| f^{\pm} - \mathcal{E}(u^{\pm}_{\varepsilon}) (u^{\pm}_{\varepsilon}) \|_{L^{3/2}(\Omega_{\pm})} \leq C \left( \| f \|_{L^{3/2}(\Omega)} + \| \nabla u^{\pm}_{\varepsilon} \|^{2}_{L^{2}(\Omega_{\pm})} \right) \, .
	\end{aligned}
\end{equation}
The pair $(u^{\pm}_{\varepsilon}, \Phi^{\pm}_{\varepsilon}) \in W^{2,3/2}(\Omega^{*}_{\pm}) \times W^{1,3/2}(\Omega^{*}_{\pm})$ is also a strong solution to the Stokes system \eqref{gammareg1}$_1$ in $\Omega^{*}_{\pm}$, with a right-hand side given by $f - \mathcal{E}(u^{\pm}_{\varepsilon}) (u^{\pm}_{\varepsilon})$. If we apply the same extension argument of the proof of \cite[Theorem 3.2]{sperone2022}, we can then invoke local regularity results for the Stokes equations (see \cite[Theorem IV.4.1]{galdi2011introduction} or the \textit{Cattabriga estimates} of \cite[Teorema, page 311]{cattabriga1961problema}) and \eqref{regsto1} to obtain the bound
\begin{equation} \label{stokesest2}
	\begin{aligned}
		& \| u^{\pm}_{\varepsilon} \|_{W^{2,3/2}(\Omega^{*}_{\pm})} + \| \Phi^{\pm}_{\varepsilon} \|_{W^{1,3/2}(\Omega^{*}_{\pm})} \\[6pt]
		& \hspace{-4mm} \leq C \left( \| f^{\pm} - \mathcal{E}(u^{\pm}_{\varepsilon}) (u^{\pm}_{\varepsilon}) \|_{L^{3/2}(\Omega_{\pm})} + \| \nabla u^{\pm}_{\varepsilon} \|_{L^{3/2}(\Omega_{\pm})} + \| \Phi^{\pm}_{\varepsilon} \|_{L^{3/2}(\Omega_{\pm})} \right) \\[6pt]
		& \hspace{-4mm} \leq C \left( 1 + \| f \|_{L^{3/2}(\Omega)} + \| \nabla u^{\pm}_{\varepsilon} \|^{2}_{L^{2}(\Omega_{\pm})} + \| \Phi^{\pm}_{\varepsilon} \|^{2}_{L^{2}(\Omega_{\pm})} \right) \qquad \forall  \varepsilon \in I_{*} \, .
	\end{aligned}
\end{equation}
Once inserted into the right-hand side of \eqref{stokesest2}, the uniform bounds in \eqref{uboundpf} entail that the sequences $(\Phi^{\pm}_{\varepsilon})_{\varepsilon \in I_{*}} \subset W^{1,3/2}(\Omega^{*}_{\pm})$ are uniformly bounded. Therefore, there exist $\widehat{\Phi}_{\pm} \in W^{1,3/2}(\Omega^{*}_{\pm})$ such that the following convergences hold as $\varepsilon \to 0^{+}$:
\begin{equation} \label{convergences3}
\Phi^{\pm}_{\varepsilon} \rightharpoonup \widehat{\Phi}_{\pm} \ \ \ \text{weakly in} \ W^{1,3/2}(\Omega^{*}_{\pm}) \qquad \text{and} \qquad \Phi^{\pm}_{\varepsilon} \to \widehat{\Phi}_{\pm} \ \ \ \text{strongly in} \ L^{1}(\partial \Omega^{*}_{\pm}) \, ,
\end{equation}
along sub-sequences that are not being relabeled, see again \cite[Theorem 6.2]{necas2011direct}. We readily notice that, since also $\Phi_{\varepsilon}^{\pm} \rightharpoonup \Phi^{\pm}$ weakly in $L^{2}(\Omega^{*}_{\pm})$ as $\varepsilon \to 0^{+}$, by uniqueness of the weak limit there must hold $\widehat{\Phi}_{\pm} = \Phi^{\pm}$ in $\Omega^{*}_{\pm}$, respectively. Moreover, in view of \eqref{gammareg1}$_2$ and the strong convergence in \eqref{convergences3}, we deduce that
\begin{equation} \label{convergences4}
\Phi^{-} = p^{-} \ \ \mbox{ almost everywhere on} \ \ \Gamma_{I} \, ; \quad \Phi^{+} = p^{+} \ \ \mbox{ almost everywhere on} \ \ \Gamma_{O} \, . 
\end{equation}
The identities in \eqref{limit1}-\eqref{limit2}-\eqref{convergences4} allow us to conclude that  $(u^{\pm}, \Phi^{\pm}) \in H^{1}(\Omega_{\pm}) \times L^{2}(\Omega_{\pm})$ are weak solutions of systems \eqref{nsstokespdfinal-}-\eqref{nsstokespdfinal+}, respectively. The improved regularity $(u^{\pm}, \Phi^{\pm}) \in H^{2}(\Omega_{\pm}) \times H^{1}(\Omega_{\pm})$ follows afterwards by repeating, again, the extension argument of \cite[Theorem 3.2]{sperone2022}. 
\par
In order to show the strong convergence in $H^{1}(\Omega_{\pm}) \times L^{2}(\Omega_{\pm})$, we start by defining the vector field $u_{*} \in H^{1}(\Omega)$ as
$$
u_{*} \doteq
\left\{ 
\begin{array}{ll}
	u^{-} & \text{in} \ \  \Omega_{-} \\[4pt]
	u^{+} & \text{in} \ \  \Omega_{+} \, ,
\end{array}
\right.
$$
(recall that $u^{-}=u^{+}=0$ on $\Sigma$) so that $u_{*}$ is divergence-free separately in $\Omega_{-}$ and $\Omega_{+}$. Fix any $\varepsilon \in I_{*}$ along the subsequences \eqref{convergences0}. Given any scalar function $\psi \in \mathcal{C}_{0}^{\infty}(\Omega)$, since $u^{-}=u^{+}=0$ on $\Sigma$, an integration by parts gives us
$$
\int_{\Omega} u_{*} \cdot \nabla \psi = \int_{\Omega_{-}} u^{-} \cdot \nabla \psi^{-} + \int_{\Omega_{+}} u^{+} \cdot \nabla \psi^{+} = \int_{\partial \Omega_{-}} \psi^{-} ( u^{-} \cdot \nu) + \int_{\partial \Omega_{+}} \psi^{+} ( u^{+} \cdot \nu) = 0 \, ,
$$
and therefore, $u_{*}$ is divergence-free in the whole $\Omega$, implying that $u_{*} \in \mathcal{V}_{\sigma}(\Omega_{\varepsilon})$. Taking $\varphi = u_{\varepsilon} - u_{*}$ as a test function in Definition \ref{weaksolution} (enforcing \eqref{refer01}) gives us
$$
\int_{\Omega_{\varepsilon}} \nabla u_{\varepsilon} \cdot \nabla (u_{\varepsilon} - u_{*}) - \int_{\Omega_{\varepsilon}} \mathcal{E}(u_{\varepsilon}) u_{\varepsilon} \cdot u_{*} + (p^{+} - p^{-}) F_{\varepsilon}  = \int_{\Omega_{\varepsilon}} f \cdot (u_{\varepsilon} - u_{*}) \, ,
$$
implying that
$$
\begin{aligned}
\| \nabla (u_{\varepsilon} - u_{*}) \|^{2}_{L^{2}(\Omega_{\varepsilon})} & = (p^{-} - p^{+}) F_{\varepsilon} + \int_{\Omega_{\varepsilon}} \nabla u_{*} \cdot \nabla (u_{*} - u_{\varepsilon}) + \int_{\Omega_{\varepsilon}} \mathcal{E}(u_{\varepsilon}) u_{\varepsilon} \cdot u_{*}  + \int_{\Omega_{\varepsilon}} f \cdot (u_{\varepsilon} - u_{*}) \\[6pt]
& = (p^{-} - p^{+}) F_{\varepsilon} + \int_{\Omega_{-}} \nabla u^{-} \cdot \nabla (u^{-} - u^{-}_{\varepsilon}) + \int_{\Omega_{-}} \mathcal{E}(u^{-}_{\varepsilon}) u^{-}_{\varepsilon} \cdot u^{-} + \int_{\Omega_{-}} f^{-} \cdot (u^{-}_{\varepsilon} - u^{-}) \\[6pt]
& \hspace{4mm} + \int_{\Omega_{+}} \nabla u^{+} \cdot \nabla (u^{+} - u^{+}_{\varepsilon}) + \int_{\Omega_{+}} \mathcal{E}(u^{+}_{\varepsilon}) u^{+}_{\varepsilon} \cdot u^{+} + \int_{\Omega_{+}} f^{+} \cdot (u^{+}_{\varepsilon} - u^{+}) \, .
\end{aligned}
$$
Taking the limit as $\varepsilon \to 0^{+}$ in this last identity, observing \eqref{tracezero}-\eqref{convergences0}, we conclude that
$$
\lim_{\varepsilon \to 0^{+}} \| \nabla (u^{-}_{\varepsilon} - u^{-}) \|_{L^{2}(\Omega_{-})} = \lim_{\varepsilon \to 0^{+}} \| \nabla (u^{+}_{\varepsilon} - u^{+}) \|_{L^{2}(\Omega_{+})} = 0 \, ,
$$
that is,
\begin{equation} \label{strong1}
u^{\pm}_{\varepsilon} \to u^{\pm} \ \ \text{strongly in} \ H^{1}(\Omega_{\pm}) \ \ \text{as} \ \ \varepsilon \to 0^{+} \, .
\end{equation}
Concerning the Bernoulli pressure, recall the sequences $(P^{\pm}_{\varepsilon})_{\varepsilon \in I_{*}} \subset L_{0}^{2}(\Omega_{\pm})$ and $(\phi^{\pm}_{\varepsilon})_{\varepsilon \in I_{*}} \subset \mathbb{R}$ introduced in \eqref{projecpres}. Likewise, define the functions $P^{\pm} \in L_{0}^{2}(\Omega_{\pm})$ by
\begin{equation} \label{projecpreslimit}
	P^{\pm} \doteq \Phi^{\pm} - \phi^{\pm} \doteq \Phi^{\pm} - \dfrac{1}{| \Omega_{\pm} |} \int_{\Omega_{\pm}} \Phi^{\pm} \qquad \text{in} \ \ \Omega_{\pm} \, .
\end{equation}
In view of the $L^{2}(\Omega_{\pm})$-weak convergence in \eqref{convergences0}, and the identities
$$
\| \Phi_{\varepsilon}^{\pm} \|^{2}_{L^{2}(\Omega_{\pm})} = \| P_{\varepsilon}^{\pm} \|^{2}_{L^{2}(\Omega_{\pm})} + | \Omega_{\pm} | | \phi_{\varepsilon}^{\pm} |^{2} \qquad \forall \varepsilon \in I_{*} \, ; \qquad \| \Phi^{\pm} \|^{2}_{L^{2}(\Omega_{\pm})} = \| P^{\pm} \|^{2}_{L^{2}(\Omega_{\pm})} + | \Omega_{\pm} | | \phi^{\pm} |^{2} \, ,
$$
it then suffices to show that
\begin{equation} \label{goal1}
\lim_{\varepsilon \to 0^{+}} \| P_{\varepsilon}^{\pm} \|_{L^{2}(\Omega_{\pm})} = \| P^{\pm} \|_{L^{2}(\Omega_{\pm})} \qquad \text{and} \qquad \lim_{\varepsilon \to 0^{+}}  \phi_{\varepsilon}^{\pm}  = \phi^{\pm}  \, .
\end{equation}
Let $(J^{\pm}_{\varepsilon})_{\varepsilon \in I_{*}} \subset H_{0}^{1}(\Omega_{\pm})$ be the sequences of vector fields satisfying \eqref{bogn0}, for every $\varepsilon \in I_{*}$ along the subsequences \eqref{convergences0}. Then, inequality \eqref{estpresind} entails that the sequences $(J^{\pm}_{\varepsilon})_{\varepsilon \in I_{*}} \subset H_{0}^{1}(\Omega_{\pm})$ must be uniformly bounded, and so we deduce the existence of vector fields $J^{\pm} \in H_{0}^{1}(\Omega_{\pm})$ such that
\begin{equation} \label{convergencesj}
	J^{\pm}_{\varepsilon} \rightharpoonup J^{\pm} \ \ \ \text{weakly in} \ H^{1}(\Omega_{\pm}) \qquad \text{and} \qquad J^{\pm}_{\varepsilon} \to J^{\pm} \ \ \ \text{strongly in} \ L^{4}(\Omega_{\pm}) \quad \text{as} \ \ \varepsilon \to 0^{+} \, ,
\end{equation}
along a (not relabeled) subsequence. Given any scalar functions $\psi^{\pm} \in \mathcal{C}_{0}^{\infty}(\Omega_{\pm})$, since $\nabla \cdot J^{\pm}_{\varepsilon} = P_{\varepsilon}^{\pm}$ in $\Omega_{\pm}$ for any $\varepsilon \in I_{*}$, an integration by parts gives us
$$
\int_{\Omega_{\pm}} J^{\pm}_{\varepsilon} \cdot \nabla \psi^{\pm} = - \int_{\Omega_{\pm}} P_{\varepsilon}^{\pm} \, \psi^{\pm} \qquad \forall \varepsilon \in I_{*} \, ,
$$
along the subsequences \eqref{convergences0}-\eqref{convergencesj}. Taking the limit in this last equality as $\varepsilon \to 0^{+}$ implies
$$
\int_{\Omega_{\pm}} J^{\pm} \cdot \nabla \psi^{\pm} = - \int_{\Omega_{\pm}} P^{\pm} \, \psi^{\pm} \qquad \forall \psi^{\pm} \in \mathcal{C}_{0}^{\infty}(\Omega_{\pm}; \mathbb{R}) \, ,
$$
that is, $\nabla \cdot J^{\pm} = P^{\pm}$ in $\Omega_{\pm}$. As in the proof of Theorem \ref{epslevelbounds} we derive the identity
\begin{equation} \label{tartar6}
\| P^{\pm}_{\varepsilon} \|^{2}_{L^{2}(\Omega_{\pm})} = \int_{\Omega_{\pm}} \nabla u^{\pm}_{\varepsilon} \cdot \nabla J^{\pm}_{\varepsilon} + \int_{\Omega_{\pm}} \mathcal{E}(u^{\pm}_{\varepsilon}) u^{\pm}_{\varepsilon} \cdot J^{\pm}_{\varepsilon} -  \int_{\Omega_{\pm}} f^{\pm} \cdot J^{\pm}_{\varepsilon} \qquad \forall \varepsilon \in I_{*} \, ,
\end{equation}
along the subsequences \eqref{convergences0}-\eqref{convergencesj}. Letting $\varepsilon \to 0^{+}$ in \eqref{tartar6}, enforcing \eqref{strong1}-\eqref{convergencesj}, we infer
\begin{equation} \label{tartar7}
\lim_{\varepsilon \to 0^{+}} \| P^{\pm}_{\varepsilon} \|^{2}_{L^{2}(\Omega_{\pm})} = \int_{\Omega_{\pm}} \nabla u^{\pm} \cdot \nabla J^{\pm} + \int_{\Omega_{\pm}} \mathcal{E}(u^{\pm}) u^{\pm} \cdot J^{\pm} -  \int_{\Omega_{\pm}} f^{\pm} \cdot J^{\pm} \, .
\end{equation}
On the other hand, observe that the identities \eqref{limit1}-\eqref{limit2} become, owing again to \eqref{split3}, 
\begin{equation} \label{limit3}
	\begin{aligned}	
		& \int_{\Omega_{-}} \nabla u^{-} \cdot \nabla \varphi + \int_{\Omega_{-}} \mathcal{E}(u^{-}) u^{-} \cdot \varphi - \int_{\Omega_{-}} P^{-} (\nabla \cdot \varphi) \\[6pt]
		& \hspace{-4mm} = \int_{\Omega_{-}} f^{-} \cdot \varphi - \left( \phi^{-} - p^{-} \right) \left( \int_{\Gamma_{I}} \varphi \cdot \widehat{k} \right) \qquad \forall \varphi \in \mathcal{V}(\Omega_{-}) \, ,
	\end{aligned}	
\end{equation}
and also
\begin{equation} \label{limit4}
	\begin{aligned}	
		& \int_{\Omega_{+}} \nabla u^{+} \cdot \nabla \varphi + \int_{\Omega_{+}} \mathcal{E}(u^{+}) u^{+} \cdot \varphi  - \int_{\Omega_{+}} P^{+} (\nabla \cdot \varphi) 	\\[6pt]
		& \hspace{-4mm}	= \int_{\Omega_{+}} f^{+} \cdot \varphi - \left( p^{+} - \phi^{+} \right) \left( \int_{\Gamma_{O}} \varphi \cdot \widehat{k} \right) \qquad \forall \varphi \in \mathcal{V}(\Omega_{+}) \, .
	\end{aligned}		
\end{equation}
Then by taking $\varphi = J^{\pm}$ as test functions in \eqref{limit3}-\eqref{limit4} (recall that $H_{0}^{1}(\Omega_{\pm}) \subset \mathcal{V}(\Omega_{\pm})$), respectively, we immediately obtain
\begin{equation} \label{tartar8}
	\| P^{\pm} \|^{2}_{L^{2}(\Omega_{\pm})} = \int_{\Omega_{\pm}} \nabla u^{\pm} \cdot \nabla J^{\pm} + \int_{\Omega_{\pm}} \mathcal{E}(u^{\pm}) u^{\pm} \cdot J^{\pm} -  \int_{\Omega_{\pm}} f^{\pm} \cdot J^{\pm} \, .
\end{equation}
Identities \eqref{tartar7}-\eqref{tartar8} secure the first limit in \eqref{goal1}, that is, 
\begin{equation} \label{strong2}
	P^{\pm}_{\varepsilon} \to P^{\pm} \ \ \text{strongly in} \ L^{2}(\Omega_{\pm}) \ \ \text{as} \ \ \varepsilon \to 0^{+} \, .
\end{equation}
Now, let us go back to the vector field $\varphi_{*} \in \mathcal{C}^{\infty}(\overline{\Omega})$ introduced in \eqref{hagen}, and formulas \eqref{split6}-\eqref{split7}, which are valid for every $\varepsilon \in I_{*}$ along the subsequences \eqref{convergences0}. Then, letting $\varepsilon \to 0^{+}$ in \eqref{split6}-\eqref{split7}, enforcing \eqref{strong1}-\eqref{strong2}, yields
\begin{equation} \label{split66}
	\begin{aligned}
	& \lim_{\varepsilon \to 0^{+}} \phi^{-}_{\varepsilon} = p^{-} + \int_{\Omega_{-}} f^{-} \cdot \varphi_{*}^{-} - \int_{\Omega_{-}} \nabla u^{-} \cdot \nabla \varphi_{*}^{-} - \int_{\Omega_{-}} \mathcal{E}(u^{-}) u^{-} \cdot \varphi_{*}^{-} + \int_{\Omega_{-}} P^{-} (\nabla \cdot \varphi_{*}^{-}) \, , \\[6pt]
	& \lim_{\varepsilon \to 0^{+}} \phi^{+}_{\varepsilon} = p^{+} - \int_{\Omega_{+}} f^{+} \cdot \varphi_{*}^{+} + \int_{\Omega_{+}} \nabla u^{+} \cdot \nabla \varphi_{*}^{+} + \int_{\Omega_{+}} \mathcal{E}(u^{+}) u^{+} \cdot \varphi_{*}^{+} - \int_{\Omega_{+}} P^{+} (\nabla \cdot \varphi_{*}^{+})  \, .
	\end{aligned}
\end{equation}
Choosing $\varphi = \varphi_{*}^{\pm}$ as test functions in \eqref{limit3}-\eqref{limit4} (recall that $\varphi_{*}^{\pm} \in \mathcal{V}(\Omega_{\pm})$), respectively, we derive
\begin{equation} \label{split67}
	\begin{aligned}
		& \phi^{-} = p^{-} + \int_{\Omega_{-}} f^{-} \cdot \varphi_{*}^{-} - \int_{\Omega_{-}} \nabla u^{-} \cdot \nabla \varphi_{*}^{-} - \int_{\Omega_{-}} \mathcal{E}(u^{-}) u^{-} \cdot \varphi_{*}^{-} + \int_{\Omega_{-}} P^{-} (\nabla \cdot \varphi_{*}^{-}) \, , \\[6pt]
		& \phi^{+} = p^{+} - \int_{\Omega_{+}} f^{+} \cdot \varphi_{*}^{+} + \int_{\Omega_{+}} \nabla u^{+} \cdot \nabla \varphi_{*}^{+} + \int_{\Omega_{+}} \mathcal{E}(u^{+}) u^{+} \cdot \varphi_{*}^{+} - \int_{\Omega_{+}} P^{+} (\nabla \cdot \varphi_{*}^{+})  \, ,
	\end{aligned}
\end{equation}
thereby ensuring the second limit in \eqref{goal1}. This concludes the proof.
\end{proof}	

In the absence of external force, the result of Theorem \ref{effectiveq1} yields a rather surprising statement:
\begin{corollary} \label{effectiveqcor}
	Let $(\Omega_{\varepsilon})_{\varepsilon \in I_{*}}$ be the family of domains \eqref{perfordomain}. Given any pressure drop $p^{\pm} \in \mathbb{R}$, let $(u_{\varepsilon},\Phi_{\varepsilon}) \in \mathcal{V}_{\sigma}(\Omega_{\varepsilon}) \times L^{2}(\Omega_{\varepsilon})$ be a weak solution of problem \eqref{nsstokespd1} (with $f=0$). Then, up to the extraction of a subsequence, the sequences of restrictions $\{(u^{\pm}_{\varepsilon}, \Phi^{\pm}_{\varepsilon}) \}_{\varepsilon \in I_{*}} \subset H^{1}(\Omega_{\pm}) \times L^{2}(\Omega_{\pm})$ converge strongly to the constant vectors $(0,p^{\pm})$ as $\varepsilon \to 0^{+}$.
\end{corollary}	
\noindent
\begin{proof}
From Theorem \ref{effectiveq1} we know that the sequences $\{(u^{\pm}_{\varepsilon}, \Phi^{\pm}_{\varepsilon}) \}_{\varepsilon \in I_{*}} \subset \mathcal{Y}_{\sigma}(\Omega_{\pm}) \times L^{2}(\Omega_{\pm})$ converge strongly to weak solutions $(u^{\pm}, \Phi^{\pm}) \in \mathcal{Y}_{\sigma}(\Omega_{\pm}) \times L^{2}(\Omega_{\pm})$ of the following boundary-value problems, in $\Omega_{\pm}$, as $\varepsilon \to 0^{+}$:
\begin{equation}\label{nsstokespdfinalcor-}
	\left\{
	\begin{aligned}
		& -\Delta u^{-} + (u^{-}\cdot\nabla)u^{-} - (\nabla u^{-})^{\top} u^{-} + \nabla \Phi^{-} = 0 \, ,\ \quad  \nabla\cdot u^{-} = 0 \ \ \mbox{ in } \ \ \Omega_{-} \, , \\[3pt]
		& u^{-} \times \nu = 0 \, , \quad \Phi^{-} = p^{-} \ \ \mbox{ on } \ \ \Gamma_{I} \, ,\\[3pt]
		& u^{-}=0 \ \ \mbox{ on } \ \ \mathcal{L}_{-} \cup \Sigma \, ,
	\end{aligned}
	\right.
\end{equation}
and
\begin{equation}\label{nsstokespdfinalcor+}
	\left\{
	\begin{aligned}
		& -\Delta u^{+} + (u^{+}\cdot\nabla)u^{+} - (\nabla u^{+})^{\top} u^{+} + \nabla \Phi^{+} = 0 \, ,\ \quad  \nabla\cdot u^{+} = 0 \ \ \mbox{ in } \ \ \Omega_{+} \, , \\[3pt]
		& u^{+} \times \nu = 0 \, , \quad \Phi^{+} = p^{+} \ \ \mbox{ on } \ \ \Gamma_{O} \, ,\\[3pt]
		& u^{+}=0 \ \ \mbox{ on } \ \ \mathcal{L}_{+} \cup \Sigma \, .
	\end{aligned}
	\right.
\end{equation}
Furthermore, $(u^{\pm}, \Phi^{\pm}) \in H^{2}(\Omega_{\pm}) \times H^{1}(\Omega_{\pm})$. The weak formulations of problems \eqref{nsstokespdfinalcor-}-\eqref{nsstokespdfinalcor+} read, respectively,
\begin{equation} \label{limitcor1}
	\int_{\Omega_{-}} \nabla u^{-} \cdot \nabla \varphi + \int_{\Omega_{-}} \mathcal{E}(u^{-}) u^{-} \cdot \varphi - p^{-} \left( \int_{\Gamma_{I}} \varphi \cdot \widehat{k} \right) - \int_{\Omega_{-}} \Phi^{-} (\nabla \cdot \varphi) = 0 \qquad \forall \varphi \in \mathcal{V}(\Omega_{-}) \, ,
\end{equation}
and also
\begin{equation} \label{limitcor2}
	\int_{\Omega_{+}} \nabla u^{+} \cdot \nabla \varphi + \int_{\Omega_{+}} \mathcal{E}(u^{+}) u^{+} \cdot \varphi + p^{+} \left( \int_{\Gamma_{O}} \varphi \cdot \widehat{k} \right) - \int_{\Omega_{+}} \Phi^{+} (\nabla \cdot \varphi) = 0 \qquad \forall \varphi \in \mathcal{V}(\Omega_{+}) \, .
\end{equation}
Recalling, from the proof of Theorem \ref{effectiveq1}, that $u^{\pm} \in \mathcal{V}(\Omega_{\pm})$ have zero flux across $\Omega_{\pm}$, it suffices to select $\varphi = u^{\pm}$ as test functions in \eqref{limitcor1}-\eqref{limitcor2}, respectively, to conclude the proof.
\end{proof}	

\begin{remark}
The result of Corollary \ref{effectiveqcor} dictates that, even under the action of an arbitrarily large prescribed pressure drop, the fluid motion becomes quiescent in the homogenization limit. 
\end{remark}

\begin{remark} \label{fremark}
	All the results presented in this manuscript hold whenever the container $\Omega$ is a smooth and distorted pipe of finite length, provided that its inlet and outlet are parallel to the $XY$-plane. To be precise, for some real numbers $0<a<b<h$, we can take
	$$
	\Omega = \left\{ (x,y,z) \in \mathbb{R}^{3} \ | \ (x,y) \in \Theta_{1} \, , \ -h < z < -a \, \right\} \cup \Omega_{c} \cup \left\{ (x,y,z) \in \mathbb{R}^{3} \ | \ (x,y) \in \Theta_{2} \, , \ b < z < h \, \right\} \, ,
	$$
	for some smooth bounded domains $\Theta_1, \Theta_2 \subset \mathbb{R}^{2}$, $\Omega_{c} \subset \mathbb{R}^{3}$ being an open, bounded, smooth and simply connected set; this geometry, depicted in Figure \ref{dom2} below, constitutes a truncation (parallel to the $XY$-plane) of the domain considered in the celebrated \textbf{Leray problem}, as described in \cite[Definition 1.1]{amick1977steady} or \cite[Chapter XIII]{galdi2011introduction}. 
	\vspace*{-2mm}
	\begin{figure}[H]
		\begin{center}
			\includegraphics[scale=0.8]{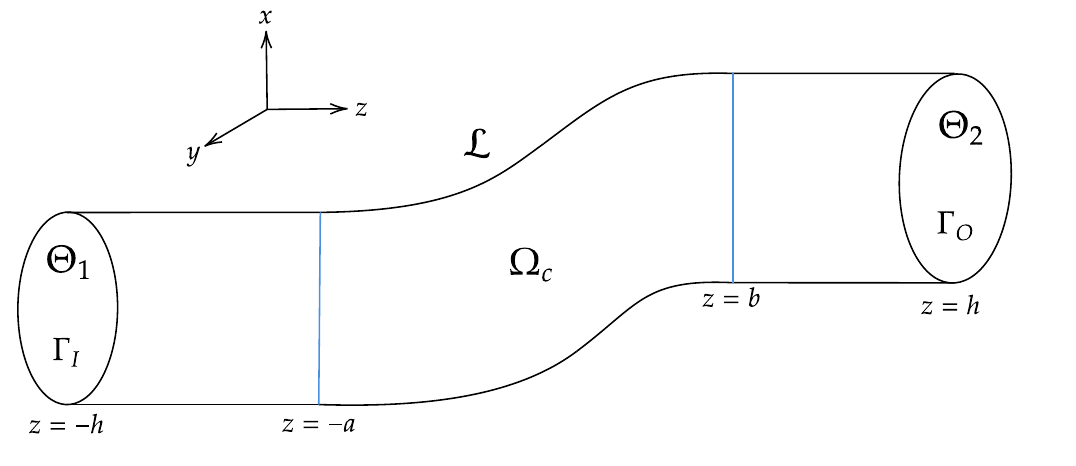}
		\end{center}
		\vspace*{-7mm}
		\caption{Representation of an admissible, smooth and distorted pipe of finite length.}\label{dom2}
	\end{figure}
	\noindent
    In this case, assuming that the cross-section $\Sigma \subset \Omega$ determinining the sieve (which is also parallel to the $XY$-plane) is located at $z=0$, the truncated Hagen-Poiseuille flow \eqref{hagen} can be defined as
	$$
	\varphi_{*}(x,y,z) = \left( \dfrac{z}{h} \right)^{2} \, V_{*}(x,y,z) \qquad \forall (x,y,z) \in \Omega \, ,
	$$
	where the vector field $V_{*} \in H^{1}(\Omega)$ is the unique weak solution to the following Stokes problem in $\Omega$:
\begin{equation} \label{stokes}
	\left\{
	\begin{aligned}
		& -\Delta V_{*} + \nabla \Pi_{*} = 0 \, , \quad  \nabla\cdot V_{*}=0 \ \ \mbox{ in } \ \ \Omega \, , \\[4pt]
		& V_{*} = V_{1} \ \ \mbox{ on } \ \ \Gamma_{I} \, , \qquad V_{*} = V_{2} \ \ \mbox{ on } \ \ \Gamma_{O} \, , \\[4pt]
		& V_{*}=0 \ \ \mbox{ on } \ \ \mathcal{L} \, . 
	\end{aligned}
	\right.
\end{equation}
In \eqref{stokes}$_2$, for $j \in \{1,2 \}$, the boundary velocity $V_{j} \in H_{0}^{1}(\Theta_{j})$ is given by
$$
V_{j}(x,y,z) \doteq \dfrac{1}{\ell_{j}} v_{j}(x,y) \widehat{k} \qquad \forall (x,y,z) \in \Theta_{j} \, ,
$$
where $v_{j} \in H^{1}_{0}(\Theta_{j}; \mathbb{R})$ is the unique weak solution of the following torsion problem in $\Theta_{j}$:
$$
-\Delta v_{j} = 1 \ \ \mbox{ in } \ \ \Theta_{j} \, , \qquad v_{j}=0 \ \ \mbox{ on } \ \ \partial \Theta_{j} \, ,
$$
and
$$
\ell_{j} \doteq \int_{\Theta_{j}} v_{j} = \int_{\Theta_{j}} | \nabla v_{j} |^{2} > 0 \, .
$$ 
\end{remark}

\par\medskip\noindent
{\bf Acknowledgements.} The Author declares that there is no conflict of interest. Data sharing not applicable to this article as no datasets were generated or analyzed during the current study.

\phantomsection
\addcontentsline{toc}{section}{References}
\bibliographystyle{abbrv}
\bibliography{references}
\vspace{5mm}

\begin{minipage}{100mm}
	Gianmarco Sperone\\
	Facultad de Matemáticas\\
	Pontificia Universidad Católica de Chile\\
	Avenida Vicuña Mackenna 4860\\
	7820436 Santiago - Chile\\
	E-mail: gianmarco.sperone@uc.cl
\end{minipage}
\end{document}